\documentclass[runningheads]{llncs}
\usepackage[utf8]{inputenc}
\usepackage[T1]{fontenc}
\usepackage{geometry,enumitem,amsmath,tikz,hyperref,color}
\usepackage[linesnumbered,boxed,ruled]{algorithm2e}
\SetAlFnt{\normalsize}

\hypersetup{
  colorlinks   = true,
  urlcolor     = blue,
  linkcolor    = blue,
  citecolor    = red
}

\def\LP{\textsc{Optimal Leaf Root}}
\def\OLR{\textsc{OLR}}
\def\dist{{\mathrm{dist}}}
\def\diam{{\mathrm{diam}}}
\def\rad{{\mathrm{rad}}}
\def\odd{{\mathrm{\pi}}}
\def\dmin{d^{\mathrm{min}}}
\def\cO{{\cal O}}
\def\cS{{\cal S}}
\def\cT{{\cal T}}

\usetikzlibrary{calc}
\tikzstyle{ivertex}=[draw,circle,inner sep=.8pt,fill=white]

\newcommand{\TPG}{ccg}
\newcommand{\TPGs}{ccgs}
\renewcommand{\oplus}{\tikz[baseline=(char.base)]{\node[ivertex] (char) {\scriptsize\normalfont 0};}}
\renewcommand{\otimes}{\tikz[baseline=(char.base)]{\node[ivertex] (char) {\scriptsize\normalfont 1};}}

\begin{document}
\title{Computing Optimal Leaf Roots of Chordal Cographs in Linear Time}
\author{Van Bang Le \and
Christian Rosenke}
\authorrunning{V.\,B.\,Le and C.\,Rosenke}

\institute{Institut f\"{u}r Informatik, Universit\"{a}t Rostock,
Rostock, Germany\\
\email{\{van-bang.le, christian.rosenke\}@uni-rostock.de}}

\maketitle

\begin{abstract}
A graph $G$ is a \emph{$k$-leaf power}, for an integer $k \geq 2$, if there is a tree $T$ with leaf set $V(G)$ such that, for all vertices $x, y \in V(G)$, the edge $xy$ exists in $G$ if and only if the distance between $x$ and $y$ in $T$ is at most $k$.
Such a tree $T$ is called a \emph{$k$-leaf root} of $G$.
The computational problem of constructing a $k$-leaf root for a given graph $G$ and an integer $k$, if any, is motivated by the challenge from computational biology to reconstruct phylogenetic trees.
For fixed $k$, Lafond [SODA 2022] recently solved this problem in polynomial time.

In this paper, we propose to study \emph{optimal leaf roots} of graphs $G$, that is, the $k$-leaf roots of $G$ with \emph{minimum} $k$ value.
Thus, all $k'$-leaf roots of $G$ satisfy $k \leq k'$.
In terms of computational biology, seeking optimal leaf roots is more justified as they yield more probable phylogenetic trees.
Lafond’s result does not imply polynomial-time computability of optimal leaf roots, because, even for optimal $k$-leaf roots, $k$ may (exponentially) depend on the size of $G$.
This paper presents a linear-time construction of optimal leaf roots for chordal cographs (also known as trivially perfect graphs).
Additionally, it highlights the importance of the parity of the parameter $k$ and provides a deeper insight into the differences between optimal $k$-leaf roots of even versus odd $k$.
\end{abstract}

\section{Introduction}
\label{sec:intro}

Leaf powers have been introduced by Nishimura, Ragde and Thilikos~\cite{NishimuraRT02} to model the phylogeny reconstruction problem from computational biology:
given a graph $G$ that represents a set of species with vertices $V(G)$ and the interspecies similarity with edges $E(G)$, how can we reconstruct an evolutionary tree $T$ with a given similarity threshold~$k$?
A \emph{$k$-leaf root} of $G$, a tree $T$ with species $V(G)$ as the leaf set and where species $x, y \in V(G)$ have distance at most $k$ in $T$ if and only if they are similar on account of $xy \in E(G)$, is considered a solution to this problem.
In case $T$ exists, the graph $G$ is called a \emph{$k$-leaf power}.
The challenge of finding a $k$-leaf root for given $G$ and $k$ has, yet, been modelled as the \emph{$k$-leaf power recognition problem}: given $G$ and $k$, decide if $G$ has a $k$-leaf root.

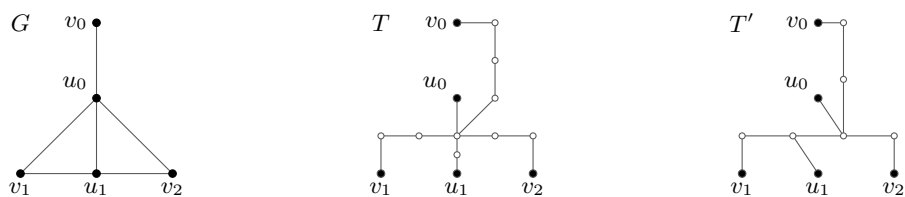
\begin{figure}\begin{center}
  \def\myscale{1}
  \tikzstyle{vertex}=[draw,circle,fill=black,inner sep=1pt]
  \tikzstyle{ivertex}=[draw,circle,inner sep=.8pt,fill=white]
  \tikzstyle{edges}=[draw=black!66, line join=round]
  \newcommand{\drawAnchors}{
    \coordinate (v1) at (0,0);
    \coordinate (u1) at (1,0);
    \coordinate (v2) at (2,0);
    \coordinate (u0) at (1,1);
    \coordinate (v0) at (1,2);
    \coordinate (l) at (v1 |- v0);
    \node[below] at (v1) {$v_1$};
    \node[below] at (u1) {$u_1$};
    \node[above left] at (u0) {$u_0$};
    \node[below] at (v2) {$v_2$};
    \node[left] at (v0) {$v_0$};
  }
  \begin{tikzpicture}[scale = \myscale, every node/.style={scale=\myscale}]
    \drawAnchors
    \node at (l) {$G$};
    \draw[edges] (v0) -- (u0) -- (v1) -- (u1) -- (u0) -- (v2) -- (u1);
    \foreach \x in {u0,v0,u1,v1,v2}{\node[vertex] at (\x) {};}
  \end{tikzpicture}
  \hspace{2cm}
  \begin{tikzpicture}[scale = \myscale, every node/.style={scale=\myscale}]
    \drawAnchors
    \node at (l) {$T$};
    \draw[edges] (v1) -- node[pos=0, vertex] {}
    ($(v1)+(0,0.5)$) -- node[pos=0, ivertex] {} node[pos=0.25, ivertex] {} node[pos=0.75, ivertex] {}
    ($(v2)+(0,0.5)$) -- 
    (v2) node[pos=0, ivertex] {} node[vertex] {};
    \draw[edges] (u1) -- node[pos=0, vertex] {} node[pos=0.5, ivertex] {}
    ($(u1)+(0,0.5)$) -- 
    (u0) node[vertex] {};
    \draw[edges] ($(u1)+(0,0.5)$) -- node[pos=0, ivertex] {}
    ($(u0)+(0.5,0)$) -- node[pos=0, ivertex] {} node[pos=0.5, ivertex] {}
    ($(v0)+(0.5,0)$) -- node[pos=0, ivertex] {}
    (v0) node[vertex] {};
  \end{tikzpicture}
  \hspace{2cm}
  \begin{tikzpicture}[scale = \myscale, every node/.style={scale=\myscale}]
    \drawAnchors
    \node at (l) {$T'$};
    \draw[edges] (v1) -- node[pos=0, vertex] {}
    ($(v1)+(0,0.5)$) -- node[pos=0, ivertex] {}
    ($(v2)+(0,0.5)$) -- node[pos=0, ivertex] {}
    (v2) node[vertex] {};
    \draw[edges] (u1) -- node[pos=0, vertex] {}
    ($(u1)+(-0.333,0.5)$) node[ivertex] {};
    \draw[edges] (u0) -- node[pos=0, vertex] {}
    ($(u0)+(0.333,-0.5)$) -- node[pos=0, ivertex] {} node[pos=0.5, ivertex] {}
    ($(v0)+(0.333,0)$) -- node[pos=0, ivertex] {}
    (v0) node[vertex] {};
  \end{tikzpicture}
  \caption{
  A graph $G$ (left), a $5$-leaf root $T$ of $G$ (middle), a $4$-leaf root $T'$ of $G$ (right). } 
\end{center}\end{figure}\label{fig:example}

\noindent
For an example, see Figure~\ref{fig:example} with the graph $G$ called \emph{dart}.
The similarities between the five species can be explained with similarity threshold $k = 5$ using the $5$-leaf root $T$ and with $k = 4$ by the $4$-leaf root $T'$, both depicted in Figure~\ref{fig:example}.

For a deeper discourse into the heavily studied field of $k$-leaf powers, the reader is kindly referred to the survey~\cite{RosenkeLB21}.
Here, we just give a short overview.

Lately, Eppstein and Havvaei~\cite{EppsteinH20} showed that $k$-leaf power recognition for graphs $G$ with $n$ vertices can be solved in $\cO(f(k, \omega) \cdot n)$ time with $f(k, \omega)$ exponential in $k$ and $\omega$, the clique number of $G$.
Quite simply put, they reduce $k$-leaf power recognition to the decision of a certain monadic second order property in a graph derived from $G$ having tree-width bounded by $k$ and $\omega$. 
Lafond's even more recent algorithm~\cite{Lafond22} solves $k$-leaf power recognition in $\cO(n^{g(k)})$ time, where $g(k)$ grows superexponentially with $k$.
It applies sophisticated dynamic programming on a tree decomposition of $G$ and exploits structural redundancies in $G$.
Observe that, for fixed $k$, the latter method runs in polynomial time.

Before these advances, $k$-leaf power recognition had only been solved for all fixed $k$ between $2$ and $6$.
The $2$-leaf powers are exactly the graphs that have just cliques as their connected components, which makes the problem trivial.
For $k=3$ (see~\cite{NishimuraRT02} and~\cite{BrandstadtL06}), $k=4$ (see~\cite{NishimuraRT02} and~\cite{BrandstadtLS08}), $k=5$ (see~\cite{ChangK07}) and $k=6$~(see~\cite{Ducoffe19}) individual algorithms have been developed, all creating a certain (tree-) decomposition of the input graph $G$ and then attempting to fit together candidate $k$-leaf roots for the components into one $k$-leaf root for $G$.

A general controversial aspect of modelling the reconstruction of a phylogenetic tree $T$ with the $k$-leaf power recognition problem is that $k$ is part of the input.
In the biological context, the value of $k$ describes an upper bound on the number of evolutionary events in $T$ that lie between two similar species $x, y$, thus, species adjacent in the given graph $G$ by an edge $xy$.
Unlike the model suggests, biologists do not always have control over the parameter $k$.
Instead, phylogenetic trees $T$ with as few as possible evolutionary events between all pairs of similar species are preferred.
That is because, in reality, a higher number of events between $x$ and $y$ makes a similarity between $x$ and $y$ less likely.
Conversely, this means that a $k$-leaf root of $G$ with a small parameter $k$ models a more probable phylogenetic tree.
This paper therefore proposes a subtle change in perspective towards considering the following optimization problem.

\bigskip\noindent
\fbox{
\begin{minipage}{.955\textwidth} 
{\LP} (\OLR)\\[.2ex]
\begin{tabular}{l p{0.85\textwidth}}
{\em Instance:}& A graph $G$.\\
{\em Output:}& An \emph{optimal leaf root} $T$ of $G$, that is,  a $\kappa$-leaf root of $G$ such that $\kappa \leq k$ for all $k$-leaf roots of $G$, or \textsc{No}, if $T$ does not exist.
\end{tabular}
\end{minipage}
}

\bigskip\noindent
Subsequently, we use $\kappa$ to indicate that the respective $\kappa$-leaf root is optimal.
{\OLR} is in a certain sense an optimization version of $k$-leaf power recognition.
The answer \textsc{No} states that the given graph $G$ is not a $k$-leaf power for any $k$ and, in particular, not for the given one.
Getting an optimal $\kappa$-leaf root of $G$ helps to decide if $G$ is a $k$-leaf power in many cases.
A difficulty is that, for all $\kappa$ and $k$ of different parity and with $2 \leq \kappa < k < 2\kappa-2$, there are $\kappa$-leaf powers that are not $k$-leaf powers~\cite{WagnerB09}.
Then, checking $\kappa \leq k$ does not decide correctly.

As for $k$-leaf power recognition, there are no known general efficient solutions for {\OLR}.
If input was restricted to $k$-leaf powers with $k \leq K$ for some fixed $K$, we could repurpose Lafond's algorithm.
Testing a given $G$ with all $2 \leq k \leq K$ would finally reveal the minimum $\kappa$ for which $G$ is $\kappa$-leaf power.
At that point, a $\kappa$-leaf root of $G$ could also be extracted from the algorithm.
But that classes of $k$-leaf powers have not been characterized well for any $k \geq 5$ makes restricting input in the proposed way difficult.
Then again, it is unknown how to decide if a given graph $G$ is a $k$-leaf power for any arbitrary $k$.
And on top of that, the minimum value $\kappa$ for which a given $G$ is a $\kappa$-leaf power, if any, may exponentially depend on the size of $G$.
This means that this brute force searching may take exponentially or even infinitely many runs of Lafond's algorithm.

It is known that, independent of $k$, all $k$-leaf powers are strongly chordal, but not vice versa.
Ptolemaic graphs are strongly chordal and a \emph{class of unbounded leaf powers}.
That is, there is no bound $\beta$ such that every Ptolemaic graph has a $k$-leaf root for some $k \leq \beta$.
Nevertheless, every Ptolemaic graph on $n$ vertices has a $2n$-leaf root~\cite{BrandstadtH08,BrandstadtHMW10}.
Later, Theorem~\ref{cor:correct_optimal} shows that, often, this is not optimal.
This paper considers a subclass of Ptolemaic graphs, the chordal cographs (also known as \emph{trivially perfect graphs}), as input to {\OLR}.
By definition, they form the intersection of the well-known \emph{chordal graphs} and the \emph{cographs}.
Accordingly, they are also characterized as the graphs without induced cycles on four vertices and without induced paths on four vertices~\cite{Golumbic78,Wolk62,Wolk65}.
As a side effect of Lemma~\ref{lem:family_roots}, this paper proves that chordal cographs are still a class of unbounded leaf powers.
This means that $k$-leaf power recognition on this class cannot be solved in polynomial time with the algorithm of Lafond or the one of Eppstein and Havvaei.
Nevertheless, the following main result of our work states that {\OLR} can be solved in linear time for chordal cographs.
\begin{theorem}\label{thm:main}
  Given a chordal cograph $G$ on $n$ vertices and $m$ edges, a (compressed) $\kappa$-leaf root of $G$ with minimum $\kappa$ can be computed in $\cO(n + m)$ time.
\end{theorem}
To the best of our knowledge, chordal cographs are, thus, the first class of unbounded leaf powers with a polynomial-time solution for {\OLR}.
The word \emph{compressed} in Theorem~\ref{thm:main} means that the $\kappa$-leaf root $T$ is returned in a denser representation, where long paths of degree two-vertices are compressed into single weighted edges.
Otherwise, the size of $T$ alone would be quadratic.

While, in general, an {\OLR}-solution does not entirely work for $k$-leaf power recognition, as elaborated above, our {\OLR}-approach can also be used for linear-time $k$-leaf power recognition on chordal cographs.
The key to this is the ability of our method to solve {\OLR} with a given parity, such that the computed $\kappa$-leaf root comes with the minimum $\kappa$ of the given parity.
Hence, if we choose the parity of the given $k$, we can tell that a given graph $G$ is a $k$-leaf power if and only if the computed $\kappa$-leaf root with $\kappa$ of the same parity as $k$ satisfies $\kappa \leq k$.

As the desired parity of $\kappa$ plays a certain role in our construction, we research this discrepancy here, and show, for certain chordal cographs, that the minimum $\kappa$ can differ up to 25 percent depending on if it is wanted odd or even.

The next section presents basic notation, definitions, and facts on trees and $k$-leaf powers used in this paper.
The optimal leaf root construction method for chordal cographs is introduced and proved correct in Section~\ref{sec:tplr}.
Section~\ref{sec:cotree} provides a respective linear-time implementation, thus, proving Theorem~\ref{thm:main}.
A deepened evaluation of the difference between chordal cographs with $\kappa$-leaf roots of minimum odd versus even $\kappa$ is carried out in the concluding Section~\ref{sec:oddeven}.

To simplify the readability of the article, all proofs have been moved to Section~\ref{sec:proofs} at the end of the paper.

\section{Preliminaries}
\label{sec:basic}

All considered graphs are finite and without multiple edges or loops.
Let $G=(V,E)$ be a graph with vertex set $V(G)=V$ and edge set $E(G)=E$. 
A \emph{universal vertex} in $G$ is one that is adjacent to all other vertices.
If all vertices of $G$ are universal than $G$ is \emph{complete}.
A vertex $x$ that is adjacent to exactly one other vertex of $G$ is called a \emph{leaf} and the edge containing $x$ is a \emph{pendant edge}.
Two adjacent vertices $x, y \in V(G)$ are \emph{true twins} if $xz \in E(G)$, if and only if $yz \in E(G)$ for all $z \in V(G) \setminus \{x, y\}$.

A graph $H$ is an \emph{induced subgraph} of $G$ if $V(H) \subseteq V(G)$ and $xy \in E(H)$ if and only if $xy \in E(G)$ for all $x,y \in V(H)$.
All subgraphs considered in this paper are induced. 
For $X \subset V(G)$, $G-X$ denotes the induced subgraph $H$ of $G$ with $V(H) = V(G) \setminus X$.
If $X$ consists of one vertex $x$ then we write $G - x$ for $G - \{x\}$.
Complete subgraphs of $G$ are called \emph{cliques}.

As usual, an \emph{$x,y$-path} in $G$ is a sequence $v_1, \dots, v_n$ of distinct vertices from $V(G)$ such that $x = v_1$, $y = v_n$ and $v_iv_{i+1} \in E(G)$ for all $i \in \{1, \dots, n-1\}$.
An $x,y$-path is called a \emph{cycle} in $G$ if $xy \in E(G)$.
The length of the $x,y$-path, respective cycle, is the number of its edges, that is, $n-1$ in the $x,y$-path and $n$ in the cycle.
If there is an $x,y$-path in $G$ for all distinct $x,y \in V(G)$ then $G$ is \emph{connected}.
Otherwise, $G$ is \emph{disconnected} and, therefore, composed of \emph{connected components} $G_1, \dots, G_n$, maximal induced subgraphs of $G$ that are connected.
A connected component is \emph{non-trivial} if it has more than one vertex and, otherwise, it is called \emph{isolated vertex}.
We call $C \subseteq V(G)$ a \emph{cut set} if $G-C$ has more connected components than $G$.
If $C$ is just a single vertex $c$ then $c$ is a \emph{cutvertex}.

Graphs $G$ and $H$ are isomorphic if a bijection $\sigma: V(G) \rightarrow V(H)$ exists with $xy \in E(G)$ if and only if $\sigma(x)\sigma(y) \in E(H)$.
If no induced subgraph of $G$ is isomorphic to a graph $H$ then $G$ is $H$-\emph{free}.
\emph{Trees} are the connected cycle-free graphs.
This means, a tree $T$ contains exactly one $x,y$-path for all $x,y \in V(T)$.

In this paper, we learn that, dependent on the given parity, the construction of an optimal leaf root differs in several details.
To avoid permanent case distinctions, we use $\odd(i)$ for the parity of an integer $i$, that is, $\odd(i)= i \mathrel{\text{mod}} 2$.

\subsection{Chordal Cographs and Cotrees}
\label{ssec:trivially_perfect}

\emph{Chordal cographs}, {\TPGs} for short, are known as the graphs that are free of the path and the cycle on $4$ vertices.
See the top row of Figure~\ref{fig:algorithm} for an example {\TPG}.
One particular {\TPG} used in this paper is the \emph{star (with $t$ leaves)}, which consists of the vertices $u, v_1, \dots, v_t$ for some $t \geq 2$ and the edges $uv_1, \dots, uv_t$.

Like all cographs, {\TPGs} can be represented with \emph{cotrees}~\cite{CorneilLB81}.
The second row of Figure~\ref{fig:algorithm} shows the cotree of the example cograph in the first row.
For every cograph $G$, the cotree $\cal T$ is a rooted tree with leaves $V(G)$ and where every internal node is labelled with $\oplus$ for \emph{disjoint union} or $\otimes$ for \emph{full join}.
In this way, the leaves define single vertex graphs and every internal node represents the cograph $G$ combining the cographs $H_1, H_2, \dots, H_n$ of its children with the respective graph operation.
More precisely, $V(G) = V(H_1) \cup V(H_2) \cup \dots \cup V(H_n)$ and $G = \oplus(H_1, H_2, \dots, H_n)$ means the disconnected cograph on vertex set $V(G)$ and edge set $E(G) = E(H_1) \cup E(H_2) \cup \dots \cup E(H_n)$ and $G = \otimes(H_1, H_2, \dots, H_n)$ means the connected cograph with vertex set $V(G)$ and edge set $E(G) = E(H_1) \cup E(H_2) \cup \dots \cup E(H_n) \cup \{xy \mid x \in V(H_i), y \in V(H_j), 1 \leq i < j \leq n\}$.
The cotree $\cal T$ is unique, can be constructed in linear time, and has the following properties:

\smallskip
\begin{itemize}[nosep]
\item
  Every internal node has at least two children.
\item
  No two internal nodes with the same label, $\oplus$ or $\otimes$, are adjacent.
\item
  The subtree ${\cal T}_X$ rooted at node $X$ is the cotree of the subgraph $G_X$ induced by the leaves of ${\cal T}_X$.
  If $X$ is labelled with $\oplus$ then $G_X$ is the disjoint union of the cographs represented by the children of $X$ and if it is labelled with $\otimes$ than $G_X$ is the full join of the children cographs.
\item
  The cotree of an $n$-vertex cograph has at most $2n-1$ nodes.
\end{itemize}

\smallskip\noindent
We mostly work with {\TPGs} without true twins, like $G$ in Figure~\ref{fig:algorithm}.
These graphs have the following properties, as observed in the upper rows of the figure.
\begin{proposition}\label{lem:trivially_perfect:cotree}
  If $G$ is a {\TPG} without true twins and $\cT$ is the cotree of $G$ then every node of $\cT$ labelled with $\otimes$ has exactly two children, one leaf and one node labelled with $\oplus$.
\end{proposition}
\begin{proposition}[Wolk~\cite{Wolk62,Wolk65}]\label{lem:wolk_plus}
  Every connected {\TPG} $G$ without true twins has a unique universal vertex $u$ and $G-u$ is disconnected (that is, $u$ is a cutvertex).
\end{proposition}

\subsection{Diameter, Radius and Center in Trees}
\label{ssec:trees}

Let $T$ be a tree.
The following notions are throughout used in the paper:

\smallskip
\begin{itemize}[nosep]
\item
The \emph{distance} between two vertices $x$ and $y$ in $T$, written $\dist_T(x,y)$, is the length of the unique $x,y$-path in $T$.
\item
The \emph{diameter} of $T$, denoted $\diam(T)$, is the maximum distance between two vertices in $T$, that is, $\diam(T) = \max\{\dist_T(x,y)\mid x, y\in V(T)\}$.
\item
A \emph{diametral path} in $T$ is a path of length $\diam(T)$.
\item
A vertex $z$ is a \emph{center vertex} of $T$ if the maximum distance between $z$ and any other vertex in $T$ is minimum, that is, all $y \in V(T)$ satisfy $\max\{\dist_T(x,z) \mid x \in V(T)\} \leq \max\{\dist_T(x,y) \mid x \in V(T)\}$.
\item
The \emph{radius} of $T$, denoted $\rad(T)$, is the maximum distance between a center $z$ and other vertices of $T$, that is, $\rad(T) = \max\{\dist_T(v,z) \mid v \in V(T)\}$. 
\end{itemize}

\smallskip\noindent
For convenience, we define $\odd(T) = \odd(\diam(T))$, the parity of the diameter of $T$.
It is well known for all trees $T$ that $\diam(T)=2\cdot\rad(T)-\odd(T)$ and that there is a single center vertex if $\odd(T) = 0$ and two adjacent centers if $\odd(T) = 1$.
Furthermore, it is obvious for any diametral path in $T$ that the end-vertices $x$ and $y$ are leaves and the center coincides with the center of $T$.
Thus, we have $\min\{\dist_T(z,x), \dist_T(z,y)\} = \rad(T)-\odd(T)$ and $\max\{\dist_T(z,x), \dist_T(z,y)\} = \rad(T)$ for any center vertex $z$ of $T$.

We call a center vertex $z$ a \emph{min-max center} of $T$ if, for all center vertices $z'$ of $T$,
$\min\{\dist_T(z,v) \mid v \text{ is a leaf of } T\} \ge \min\{\dist_T(z',v) \mid v \text{ is a leaf of } T\}$.
Thus, a min-max center maximizes the distance to the closest leaf of $T$. 
For a min-max center $z$, the \emph{leaf distance} is $\dmin_T = \min\{\dist_T(z,v) \mid v \text{ is a leaf of } T\}$.
The paper also needs the following technical lemma:
\begin{lemma}\label{lem:technical}
If $T_1, \dots, T_s$ are $s \ge 2$ trees with $\diam(T_1) \ge \dots \ge \diam(T_s)$ then
\begin{enumerate}[label=\textup{(\roman*)},leftmargin=*,nosep]
\item\label{lem:technical:i}
$\rad(T_1)\ge \rad(T_2)\ge \cdots \ge \rad(T_s)$,
\item\label{lem:technical:ii}
for all $1\le i<j\le s$, if $\rad(T_i)=\rad(T_j)$ then $\odd(T_i)\le \odd(T_j)$, and
\item\label{lem:technical:iii}
for all $1\le i<j \le s$, $\rad(T_i)-\odd(T_i)\ge \rad(T_j)-\odd(T_j)$. 
\end{enumerate}
\end{lemma}

\subsection{Leaf Powers, Leaf Roots and their Basic Properties}
\label{ssec:lp}

Let $k \ge 2$ be an integer.
A graph $G$ is a \emph{$k$-leaf power} if a \emph{$k$-leaf root} $T$ of $G$ exists, a tree with leaves $V(G)$ such that $xy$ is an edge in $G$ if and only if $\dist_T(x,y) \leq k$.
The example $G$ in Figure~\ref{fig:algorithm} therefore is an $11$-leaf power because of the $11$-leaf root $T$ of $G$ in the third row and also a $12$-leaf power by the $12$-leaf root $T'$ in the bottom row.
Note that Figure~\ref{fig:algorithm} shows \emph{compressed} illustrations of $T$ and $T'$ where some long paths of vertices with degree two are depicted by single weighted edges.
It is well-known that

\smallskip
\begin{itemize}[nosep]
\item a complete graph is a $k$-leaf power for all $k \geq 2$,
\item a graph is a $k$-leaf power if and only if all of its connected components are $k$-leaf powers, and
\item if $x,y$ are true twins in $G$ then $G$ is a $k$-leaf power if and only if $G-x$ is a $k$-leaf power.
\end{itemize}

\smallskip\noindent
Note for the last fact that Lemma~7.3 and Corollary~7.4 in~\cite{McConnell03} imply the possibility to identify and remove all true twins from a graph in linear time.
So, in the remainder of the paper, we smoothly focus on graphs without true twins.

Since the concept of $k$-leaf powers is slightly different for odd and even $k$, we formalize this discrepancy as follows:
We say that a $k$-leaf root is of \emph{even}, respectively \emph{odd parity}, if $k$ is even, respectively odd.
A $k$-leaf root $T$ of $G$ is an \emph{optimal even}, respectively \emph{optimal odd} leaf root if $k$ is even, respectively odd, and all $k'$-leaf roots of $G$ with $k'$ of the same parity as $k$ satisfy $k \le k'$.
Finally, a $k$-leaf root $T$ of $G$ is (just) \emph{optimal} if $k \le k'$ for all $k'$-leaf roots of $G$ (independent of parity).
See the third row of Figure~\ref{fig:algorithm} for an optimal odd leaf root $T$ of the example graph $G$ in the same figure and see the bottom row for an optimal even leaf root $T'$ of $G$.
Since $T$ is an $11$-leaf root and $T'$ a $12$-leaf root, it follows that $T$ is an optimal leaf root of $G$.

We conclude this section by establishing a few properties for leaf roots as considered in this paper.
The first one concerns a bound on the distance between center and leaves in $T$ in case $G$ is connected.
\begin{lemma}\label{lem:dmin}
  Every $k$-leaf root $T$ of a connected graph satisfies $\dmin_T \le \tfrac{k}{2}$. 
\end{lemma}
See Figure~\ref{fig:algorithm} with the $11$-leaf root $T$ of $G$ having $\dmin_T = \dist_T(u_0, z_0) = 2 \leq \tfrac{11}{2}$ and the $12$-leaf root $T'$ with $\dmin_{T'} = \dist_T(u_0, z_0) = 2 \leq \tfrac{12}{2}$. 
Secondly, if $G$ has a universal vertex $u$ then the distance between $u$ and the center in $T$ cannot exceed the difference between $k$ and the radius of $T$.
\begin{lemma}\label{lem:A} 
If $G$ is a non-complete graph with a universal vertex $u$ and $T$ is a $k$-leaf root of $G$ then $\dist_T(u,z) \le k-\rad(T)+\odd(T)$ for all center vertices $z$ of $T$.
If $z_1 \not= z_2$ are the center vertices of $T$, then $\dist_T(u,z_1) \leq k-\rad(T)$ or $\dist_T(u,z_2) \leq k-\rad(T)$.
\end{lemma}
For an illustration, see Figure~\ref{fig:algorithm}, where the distance of $u_0$ and the farthest center vertex $z_0$ of $T$ satisfies $\dist_T(u_0,z_0) = 2 \leq 11 - 10 + 1 = k - \rad(T) + \odd(T)$ and, in $T'$, $\dist_{T'}(u_0,z_0) = 2 \leq 12 - 11 + 1 = k' - \rad(T') + \odd(T')$. 
Lemma~\ref{lem:A} implies upper bounds on radius and diameter of $T$.
\begin{corollary}\label{cor:A}  
  If $G$ is a graph with a universal vertex and $T$ is a $k$-leaf root of $G$ then $\rad(T)\le k-1$ and, in particular, $\diam(T)\le 2k-2$.
\end{corollary}
As a matter of fact, $k$-leaf roots tend to contain long paths of vertices with degree two.
It is reasonable to \emph{compress} such a path $P = v_0, \dots, v_n$ into a single \emph{weighted} edge $v_0(n)v_n$.
Clearly, weighted edges $v_0(n)v_n$ add their weight $n$ to distances in $T$ and, so, $\dist_T(v_0, v_n) = n$.

\section{Optimal Leaf Root Construction for CCGs}
\label{sec:tplr}

Aim of this section is the development of an optimal leaf-root construction approach for {\TPGs} $G$.
In very simple terms, we describe a divide and conquer method that splits $G$ into smaller {\TPGs} $G_1, G_2, \dots$, recursively obtains their optimal leaf roots $T_1, T_2, \dots$, and then extends them into an optimal leaf root for $G$.

We start with introducing two basic leaf root operations and analyze their properties.
The first operation, the \emph{extension} of trees, is used to level the recursively found leaf roots $T_1, T_2, \dots$ on the same $k$, which is essential for the subsequent composition into one $k$-leaf root.
If $T$ is a tree and $\delta \geq 0$ an integer then $T' = \eta(T, \delta)$ is the obtained from $T$ by subdividing every pendant edge $\delta$ times, that is, replacing the edge with a new path of length $\delta+1$ (hence, of $\delta+1$ edges).
The following property of this operation is well-known:
\begin{lemma}\label{lem:subdividing}
If $T$ is a $k$-leaf root of a graph $G$ and $\delta \ge 0$ an integer then $T' = \eta(T, \delta)$ is a $(k+2\delta)$-leaf root of $G$ with same center, same min-max center vertices, and
$\diam(T') = \diam(T)+2\delta$, $\rad(T') = \rad(T)+\delta$, $\dmin_{T'}=\dmin_{T}+\delta$.
\end{lemma}
The second operation merges the individual $k$-leaf roots for the connected components of a graph $G$ into one $k$-leaf root $T$ for the entire $G$.
The goal here is to minimize the diameter of $T$, which, in turn, allows making optimizations to the value of $k$.
Assume that $G$ has $s \ge 0$ non-trivial connected components $G_1, \dots, G_s$ and $t \ge 0$ isolated vertices $v_1, \dots, v_t$ such that $s+t \ge 2$ and let $T_1, \dots, T_s$ be $k$-leaf roots for $G_1, \dots, G_s$ with min-max center vertices $z_1, \dots, z_s$.
If $s > 0$, we define the \emph{critical index} $m$ as the smallest element of $\{1, \dots, s\}$ with $\dmin_{T_m} = \min\{\dmin_{T_i} \mid 1 \leq i \leq s\}$ and call $T_m$ the \emph{critical root}.
Then, the \emph{merging} $\mu(k, T_1, \dots, T_s, v_1, \dots, v_t)$ results in the tree $T$ produced by the following steps:

\smallskip
\begin{enumerate}[nosep]
  \item Create a new vertex $c$.
  \item If $s > 0$ then connect $c$ and the center $z_m$ of the critical root by a path of length $\tfrac{k+\odd(k)}{2} - \dmin_{T_m}$.
  If $\odd(k) = 0$ and $\dmin_{T_m} = \tfrac{1}{2}k$ then this means to identify the vertices $c$ and $z_m$.
  \item For all $i \in \{1, \dots, m-1, m+1, \dots, s\}$, connect $c$ and $z_i$ by a path of length $\tfrac{k-\odd(k)}{2}+1-\dmin_{T_i}$.
  \item For all $j \in \{1, \dots, t\}$, connect $c$ and $v_j$ by a path of length $\tfrac{k-\odd(k)}{2}+1$.
\end{enumerate}

\smallskip\noindent
Notice that the $\mu$-operation is sensitive with respect to the parity $\odd(k)$.
For one thing, this is necessary to guarantee that all added paths are of integer length, which is done by in- or decreasing odd $k$.
As a side note, we point out that the lengths of added paths are also non-negative by Lemma~\ref{lem:dmin}, which makes the $\mu$-operation well-defined.
On the other hand, the result is that merging works slightly different for odd and even $k$.
For odd $k$, all trees $T_1, \dots, T_s$, including the critical one, are essentially added in the same way by our construction method.
This is because, for odd $k$, $\tfrac{k+\odd(k)}{2} - \dmin_{T_m} = \tfrac{k-\odd(k)}{2}+1-\dmin_{T_m}$.
The special treatment of the critical root, thus, has an effect only if $k$ is even.
Specifically in that case, we can sometimes save one in the diameter of $T$, if we put the critical root closer to the center of $T$ than the rest.
The reason that this optimization works is that, usually, the critical root has the largest diameter.
\begin{lemma}\label{lem:disconnected}
  Let $G$ be a graph and $k \geq 2$ an integer.
  If $G$ is disconnected with $s \ge 0$ non-trivial connected components $G_1, \dots, G_s$ and $t \ge 0$ isolated vertices $v_1, \dots, v_t$ such that $s+t \ge 2$ and if $T_1, \dots, T_s$ are $k$-leaf roots for $G_1, \dots, G_s$ then $T = \mu(k, T_1, \dots, T_s, v_1, \dots, v_t)$ is a $k$-leaf root of $G$.
\end{lemma}
The two operations above simplify the description of the following leaf root construction algorithm for {\TPGs} since they hide away many of the technical details.
Foundation of the proposed recursive approach is that (i) induced subgraphs of {\TPGs} are {\TPGs} and (ii) every connected {\TPG} without true twins has a unique universal cut vertex (see Proposition~\ref{lem:wolk_plus}).
Also, recall from Section~\ref{ssec:lp} that true twins in graphs can be removed in linear time and, thus, be safely ignored.
Therefore, we define for all {\TPGs} $G$ without true twins and a given parity $p \in \{0,1\}$ the result of the \emph{root operation} $\rho(G, p)$ as the tree $T$ and the number $k$ produced by the following (recursive) procedure:

\smallskip
\begin{enumerate}[label=\roman*.,nosep]
  \item\label{const:base}\textbf{If $G$ is a star then}
  let $u$ be the central vertex and $v_1, \dots, v_t$ the leaves of $G$ (with $t \geq 2$ because $G$ does not have true twins)
  \begin{enumerate}[label=\arabic*.,nosep]
    \item\label{const:base:odd} 
    If $p=1$ (for odd) then let $T' = \eta(G, 1)$, obtain $T$ by attaching a new leaf to $u$ in $T'$, and return $(T, 3)$.
    \item\label{const:base:even2} 
    If $p=0$ (for even) and $t=2$ then return $(T, 4)$ with $T$ obtained from a single vertex $v$ by attaching the leaves $u, v_1$ and $v_2$ to $v$ with paths of lengths one, two and three, respectively.
    \item\label{const:base:evenN} 
    If $p=0$ and $t>2$ then let $T' = \eta(G, 2)$, obtain $T$ by attaching a new leaf to $u$ in $T'$, and return $(T, 4)$.
  \end{enumerate}
  \item\label{const:case:connected}\textbf{else if $G$ is a connected graph then}
  let $u$ be the universal cut vertex of $G$ (by Proposition~\ref{lem:wolk_plus}) and let $G_1, \dots, G_s$ be the $s \ge 1$ non-trivial connected components and $v_1, \dots, v_t$ the $t \ge 0$ isolated vertices of $G-u$.
  \begin{enumerate}[label=\arabic*.,nosep]
    \item\label{const:case:connected:recursion}
    Recursively find $(T_1, k_1) = \rho(G_1, p), \dots, (T_s, k_s) = \rho(G_s, p)$.
    \item\label{const:case:connected:k} 
    If $s = 1$ then let $k = k_1 + 2(1 - \odd(T_1))$.
    Otherwise, let
    
    \smallskip
    \centerline{$\begin{aligned}
      k_a &= \max\{k_1, \dots, k_s\},\\
      k_b &= \max\{k_i \mid 1 \leq i \leq s, i \not= a\} \text{ and, if } s > 2 \text{ let}\\
      k_c &= \max\{k_i \mid 1 \leq i \leq s, i \not= a, i \not= b\}.
    \end{aligned}$}

    \smallskip
    If $p=1$ (for odd) then let $k = k_a + k_b - 1 - 2\cdot\odd(T_a)\cdot\odd(T_b)$ and, otherwise,
    \[k = \begin{cases}
      k_a + k_b - 2\cdot(\odd(T_a)+\odd(T_b)-\odd(T_a)\cdot\odd(T_b)), & \text{if } s=2 \text{ or } k_a > k_c\\
      k_a + k_b - 2\cdot\odd(T_a)\cdot\odd(T_b), &\text{otherwise.}
    \end{cases}\]
    \item\label{const:case:connected:widen}
    Get the extended leaf root $T'_i = \eta\left(T_i, \tfrac{k-k_i}{2}\right)$ for all $i \in \{1, \dots, s\}$ and let $T' = \mu(k, T'_1, \dots, T'_s, v_1, \dots, v_t)$.
    \item\label{const:case:connected:assembly} 
    Return $(T,k)$ with $T$ obtained from $T'$ by attaching the leaf $u$ to a center vertex of $T'$.
  \end{enumerate}
  \item\label{const:case:disconnected}\textbf{else $G$ is a disconnected graph and then}
  let $G_1, \dots, G_s$ be the $s \ge 0$ non-trivial connected components and $v_1, \dots, v_t$ the $t \ge 0$ isolated vertices of $G$.
  \begin{enumerate}[label=\arabic*.,nosep]
    \item\label{const:case:disconnected:recursion} 
    Recursively find $(T_1, k_1) = \rho(G_1, p), \dots, (T_s, k_s) = \rho(G_s, p)$.
    \item\label{const:case:disconnected:widen}
    Let $k = \max\{k_1, \dots, k_s, p+2\}$ and let $T'_i = \eta\left(T_i, \tfrac{k-k_i}{2}\right)$ for all $i \in \{1, \dots, s\}$.
    \item\label{const:case:disconnected:assembly} 
    Return $(T, k)$ with $T = \mu(k, T'_1, \dots, T'_s, v_1, \dots, v_t)$.
  \end{enumerate}
\end{enumerate}

\begin{figure}\begin{center}
  \def\myscale{0.66}
  \tikzstyle{vertex}=[draw,circle,fill=black,inner sep=1pt]
  \tikzstyle{ivertex}=[draw,circle,inner sep=.8pt,fill=white]
  \tikzstyle{edges}=[draw=black!66, line join=round]
  \tikzstyle{subgraph}=[fill=black,rounded corners, fill opacity=0.10] 
  \newcommand{\drawAnchorsAndBackground}{
    \coordinate (gn1) at (1,0);
    \coordinate (gn0) at (0,-1.5);
    \coordinate (u0) at (0,-0.5);
    \coordinate (v0) at (0,-2.5);
    \coordinate (g1n1) at (-4,-2.5);
    \coordinate (g1n0) at (-5,-3.5);
    \coordinate (u1) at (-5,-3);
    \coordinate (g2n1) at (4,-2.5);
    \coordinate (g2n0) at (5,-3.5);
    \coordinate (u2) at (5,-3);
    \coordinate (g3n1) at (-7,-4);
    \coordinate (g3n0) at (-8,-5);
    \coordinate (u3) at (-8,-4.5);
    \coordinate (v1) at (-7,-5.5);
    \coordinate (g4n1) at (-4,-4.5);
    \coordinate (g4n0) at (-5,-5.5);
    \coordinate (u4) at (-5,-5);
    \coordinate (v2) at (-6,-6);
    \coordinate (v3) at (-5,-6);
    \coordinate (v4) at (-4,-6);
    \coordinate (g5n1) at (-1,-4.5);
    \coordinate (g5n0) at (-2,-5.5);
    \coordinate (u5) at (-2,-5);
    \coordinate (v5) at (-3,-6);
    \coordinate (v6) at (-1,-6);
    \coordinate (g6n1) at (3,-4.5);
    \coordinate (g6n0) at (2,-5.5);
    \coordinate (u6) at (2,-5);
    \coordinate (v7) at (1,-6);
    \coordinate (v8) at (3,-6);
    \coordinate (g7n1) at (6,-4.5);
    \coordinate (g7n0) at (5,-5.5);
    \coordinate (u7) at (5,-5);
    \coordinate (v9) at (4,-6);
    \coordinate (v10) at (6,-6);
    \coordinate (g8n1) at (9,-4.5);
    \coordinate (g8n0) at (8,-5.5);
    \coordinate (u8) at (8,-5);
    \coordinate (v11) at (7,-6);
    \coordinate (v12) at (9,-6);
    \coordinate (g9n1) at (-7,-6);
    \coordinate (g9n0) at (-8,-7);
    \coordinate (u9) at (-8,-6.5);
    \coordinate (v13) at (-9,-7.5);
    \coordinate (v14) at (-7,-7.5);
    \fill[subgraph] ($(v13)+(-0.8,-0.8)$) rectangle ($(g1n1 -| v6)+(0.7,0.7)$);
    \fill[subgraph] ($(v7)+(-0.6,-1.6)$) rectangle ($(g2n1 -| v12)+(0.6,0.6)$);
    \fill[subgraph] ($(v13)+(-0.6,-0.6)$) rectangle ($(g3n1)+(0.5,0.5)$);
    \fill[subgraph] ($(v2)+(-0.4,-1.4)$) rectangle ($(g4n1)+(0.4,0.4)$);
    \fill[subgraph] ($(v5)+(-0.4,-1.4)$) rectangle ($(g5n1)+(0.4,0.4)$);
    \fill[subgraph] ($(v7)+(-0.4,-1.4)$) rectangle ($(g6n1)+(0.4,0.4)$);
    \fill[subgraph] ($(v9)+(-0.4,-1.4)$) rectangle ($(g7n1)+(0.4,0.4)$);
    \fill[subgraph] ($(v11)+(-0.4,-1.4)$) rectangle ($(g8n1)+(0.4,0.4)$);
    \fill[subgraph] ($(v13)+(-0.4,-0.4)$) rectangle ($(g9n1)+(0.3,0.3)$);
    \coordinate (l0) at (v13 |- u0);
    \coordinate (l1) at (v13 |- u1);
    \coordinate (l2) at (v12 |- u2);
    \coordinate (l3) at (v13 |- u3);
    \coordinate (l4) at ($(v2) + (0,-1)$);
    \coordinate (l5) at ($(v5) + (0,-1)$);
    \coordinate (l6) at ($(v7) + (0,-1)$);
    \coordinate (l7) at ($(v9) + (0,-1)$);
    \coordinate (l8) at ($(v11) + (0,-1)$);
    \coordinate (l9) at (v13 |- u9);
  }
  %
  \begin{tikzpicture}[scale = \myscale, every node/.style={scale=\myscale}]
    \drawAnchorsAndBackground
    \draw[edges] (u3) -- (u1) -- (u4) -- (u0) -- (u5) -- (u1) -- (u0) -- (v6) -- (u5) -- (v5) -- (u1) -- (v4) -- (u4) -- (v2) -- (u1) -- (v1) -- (u3) -- (v13) -- (u9) -- (v14) -- (u3) -- (u9);
    \draw[edges] (u4) -- (v3) -- (u0) -- (v4);
    \draw[edges] (v1) -- (u0) -- (v11);
    \draw[edges] (v4) -- (u0) -- (v5);
    \draw[edges] (u0) -- (v0);
    \draw[edges] (v9) -- (u0) -- (v8);
    \draw[edges] (u8) -- (u2) -- (u7) -- (u0) -- (u6) -- (u2) -- (u0) -- (v7) -- (u6) -- (v8) -- (u2) -- (v9) -- (u7) -- (v10) -- (u2) -- (v11) -- (u8) -- (v12);
    \draw[edges] (u0) to[bend angle=10, bend right] (u3);
    \draw[edges] (u0) .. controls ($(u0)+(-1,0)$) and ($(v13)+(0,7)$) .. (v13);
    \draw[edges] (u0) .. controls ($(v1)+(0,0.5)$) .. (u9);
    \draw[edges] (u0) .. controls ($(v2)+(0,1)$) .. (v14);
    \draw[edges] (u0) .. controls ($(u4)+(-0.2,0.2)$) .. (v2);
    \draw[edges] (u0) .. controls ($(u7)+(0.2,0.2)$) .. (v10);
    \draw[edges] (u0) to[bend angle=15, bend left] (u8);
    \draw[edges] (u0) to[bend angle=30, bend left] (v12);
    \draw[edges] (u1) to[bend angle=10, bend right] (v3);
    \draw[edges] (u1) to[bend angle=20, bend left] (v6);
    \draw[edges] (u1) to[bend angle=15, bend right] (v13);
    \draw[edges] (u1) to[bend angle=10, bend right] (u9);
    \draw[edges] (u1) to[bend angle=5, bend right] (v14);
    \draw[edges] (u2) to[bend angle=20, bend left] (v12);
    \draw[edges] (u2) to[bend angle=20, bend right] (v7);
    \foreach \u in {0,1,...,9}{\node[vertex] at (u\u) {};}
    \foreach \u in {0,1,2}{\node[above] at (u\u) {$u_{\u}$};}
    \node[above=0.075cm] at (u3) {$u_{3}$};
    \node[right=0.1cm] at (u4) {$u_{4}$};
    \foreach \u in {5,6,7,8,9}{\node[below=0.1cm] at (u\u) {$u_{\u}$};}
    \foreach \v in {0,1,...,14}{
      \node[vertex] at (v\v) {};
      \node[below] at (v\v) {$v_{\v}$};
    }
    \node at (l0) {$G$};
    \node at (l1) {$G_1$};
    \node at (l2) {$G_2$};
    \node at (l3) {$G_3$};
    \node at (l4) {$G_4$};
    \node at (l5) {$G_5$};
    \node at (l6) {$G_6$};
    \node at (l7) {$G_7$};
    \node at (l8) {$G_8$};
    \node at (l9) {$G_9$};
  \end{tikzpicture}
  
  \medskip
  \begin{tikzpicture}[scale = \myscale, every node/.style={scale=\myscale}]
    \drawAnchorsAndBackground
    \draw[edges] (u0) -- (gn1) -- (gn0) -- (g1n1);
    \draw[edges] (v0) -- (gn0) -- (g2n1);
    \draw[edges] (u1) -- (g1n1) -- (g1n0) -- (g3n1);
    \draw[edges] (u2) -- (g2n1) -- (g2n0) -- (g6n1);
    \draw[edges] (u3) -- (g3n1) -- (g3n0) -- (v1);
    \draw[edges] (g3n0) -- (g9n1);
    \draw[edges] (g4n1) -- (g1n0) -- (g5n1);
    \draw[edges] (g7n1) -- (g2n0) -- (g8n1);
    \draw[edges] (u4) -- (g4n1) -- (g4n0) -- (v2);
    \draw[edges] (v3) -- (g4n0) -- (v4);
    \draw[edges] (u5) -- (g5n1) -- (g5n0) -- (v5);
    \draw[edges] (v6) -- (g5n0);
    \draw[edges] (u6) -- (g6n1) -- (g6n0) -- (v7);
    \draw[edges] (v8) -- (g6n0);
    \draw[edges] (u7) -- (g7n1) -- (g7n0) -- (v9);
    \draw[edges] (v10) -- (g7n0);
    \draw[edges] (u8) -- (g8n1) -- (g8n0) -- (v11);
    \draw[edges] (v12) -- (g8n0);
    \draw[edges] (u9) -- (g9n1) -- (g9n0) -- (v13);
    \draw[edges] (v14) -- (g9n0);
  
    \node[ivertex] at (gn0) {0};
    \node[ivertex] at (gn1) {1};
    \foreach \g in {1,...,9}{
      \node[ivertex] at (g\g n0) {0};
      \node[ivertex] at (g\g n1) {1};
    }
    \foreach \u in {0,1,...,9}{\node[vertex] at (u\u) {};}
    \foreach \u in {0,1,3,4,5,6,7,8,9}{\node[left] at (u\u) {$u_{\u}$};}
    \node[right] at (u2) {$u_{2}$};
    \foreach \v in {0,1,...,14}{\node[vertex] at (v\v) {};}
    \foreach \v in {0,2,3,4,5,6,7,8,9,10,11,12,13,14}{\node[below] at (v\v) {$v_{\v}$};}
    \node[right] at (v1) {$v_{1}$};
    \node at (l0) {$\cT$};
    \node at (l1) {$\cT_1$};
    \node at (l2) {$\cT_2$};
    \node at (l3) {$\cT_3$};
    \node at (l4) {$\cT_4$};
    \node at (l5) {$\cT_5$};
    \node at (l6) {$\cT_6$};
    \node at (l7) {$\cT_7$};
    \node at (l8) {$\cT_8$};
    \node at (l9) {$\cT_9$};
  \end{tikzpicture}
  
  \medskip
  \newcommand{\subdividedEdge}[4]{
    \draw[edges] (#1) #4 (#2);
  }
  \newcommand{\leafEdge}[5]{
    \draw[edges] (#1) #5 node[pos=0.5,#4]{\tiny $#3$} (#2);
  }
  \begin{tikzpicture}[scale = \myscale, every node/.style={scale=\myscale}]
    \drawAnchorsAndBackground
    \subdividedEdge{$(gn0)+(-0.5,0)$}{$(gn0)+(0.5,0)$}{1}{--}
    \subdividedEdge{$(gn0)+(-0.5,0)$}{g1n0}{1}{--}
    \subdividedEdge{$(gn0)+(0.5,0)$}{g2n0}{1}{-- node[ivertex, pos=0.5] {}}
    \leafEdge{$(gn0)+(-0.5,0)$}{u0}{}{right=-0.05cm}{--}
    \leafEdge{$(gn0)+(0.5,0)$}{v0}{6}{below right=-0.05cm}{--}
    \node[ivertex] at ($(gn0)+(-0.5,0)$) {};
    \node[ivertex] at ($(gn0)+(0.5,0)$) {};
    \node[above] at ($(gn0)+(0.5,0)$) {$z_0$};
    \subdividedEdge{g1n0}{g3n0}{1}{--}
    \subdividedEdge{g1n0}{g4n0}{1}{to[bend angle=25, bend left]}
    \subdividedEdge{g1n0}{g5n0}{1}{--}
    \leafEdge{g3n0}{u1}{3}{above}{to[bend angle=25, bend left]}
    \node[ivertex] at (g1n0) {};
    \node[below=0.1cm] at ($(g1n0)+(-0.1,0)$) {$z_1$};
    \subdividedEdge{g2n0}{g6n0}{1}{--}
    \subdividedEdge{g2n0}{g7n0}{1}{to[bend angle=25, bend left]}
    \subdividedEdge{g2n0}{g8n0}{1}{--}
    \leafEdge{g2n0}{u2}{4}{right=-0.05cm}{--}
    \node[ivertex] at (g2n0) {};
    \node[below=0.1cm] at ($(g2n0)+(-0.1,0)$) {$z_2$};
    \leafEdge{g3n0}{v1}{6}{above right=-0.05cm}{--}
    \subdividedEdge{g3n0}{g9n0}{6}{to[bend angle=25, bend right]}
    \leafEdge{g9n0}{u3}{4}{left}{to[bend angle=40, bend left]}
    \node[ivertex] at (g3n0) {};
    \node[below=0.1cm] at ($(g3n0)+(0.1,0)$) {$z_3$};
    \leafEdge{g4n0}{u4}{5}{left=-0.05cm}{--}
    \leafEdge{g4n0}{v2}{6}{above}{--}
    \leafEdge{g4n0}{v3}{6}{right=-0.05cm}{--}
    \leafEdge{g4n0}{v4}{6}{above}{--}
    \node[ivertex] at (g4n0) {};
    \node[below=0.1cm] at ($(g4n0)+(-0.2,0)$) {$z_4$};
    \leafEdge{g5n0}{u5}{5}{right=-0.05cm}{--}
    \leafEdge{g5n0}{v5}{6}{above}{--}
    \leafEdge{g5n0}{v6}{6}{above}{--}
    \node[ivertex] at (g5n0) {};
    \node[below=0.1cm] at (g5n0) {$z_5$};
    \leafEdge{g6n0}{u6}{5}{left=-0.05cm}{--}
    \leafEdge{g6n0}{v7}{6}{above}{--}
    \leafEdge{g6n0}{v8}{6}{above}{--}
    \node[ivertex] at (g6n0) {};
    \node[below=0.1cm] at (g6n0) {$z_6$};
    \leafEdge{g7n0}{u7}{5}{left=-0.05cm}{--}
    \leafEdge{g7n0}{v9}{6}{above}{--}
    \leafEdge{g7n0}{v10}{6}{above}{--}
    \node[ivertex] at (g7n0) {};
    \node[below=0.1cm] at (g7n0) {$z_7$};
    \leafEdge{g8n0}{u8}{5}{right=-0.05cm}{--}
    \leafEdge{g8n0}{v11}{6}{above}{--}
    \leafEdge{g8n0}{v12}{6}{above}{--}
    \node[ivertex] at (g8n0) {};
    \node[below=0.1cm] at (g8n0) {$z_8$};
    \leafEdge{g9n0}{u9}{5}{right=-0.05cm}{--}
    \leafEdge{g9n0}{v13}{6}{above}{--}
    \leafEdge{g9n0}{v14}{6}{above}{--}
    \node[ivertex] at (g9n0) {};
    \node[below=0.1cm] at (g9n0) {$z_9$};
    \foreach \u in {0,1,...,9}{
      \node[vertex] at (u\u) {};
      \node[above] at (u\u) {$u_{\u}$};
    }
    \foreach \v in {0,1,...,14}{
      \node[vertex] at (v\v) {};
      \node[below] at (v\v) {$v_{\v}$};
    }
    \node at (l0) {$T$};
    \node at (l1) {$T_1$};
    \node at (l2) {$T_2$};
    \node at (l3) {$T_3$};
    \node at (l4) {$T_4$};
    \node at (l5) {$T_5$};
    \node at (l6) {$T_6$};
    \node at (l7) {$T_7$};
    \node at (l8) {$T_8$};
    \node at (l9) {$T_9$};
  \end{tikzpicture}
  
  \medskip
  \begin{tikzpicture}[scale = \myscale, every node/.style={scale=\myscale}]
    \drawAnchorsAndBackground
    \subdividedEdge{$(gn0)+(-0.5,0)$}{g1n1}{1}{-- node[ivertex, pos=0.5]{}}
    \subdividedEdge{$(gn0)+(-0.5,0)$}{$(gn0)+(0.5,0)$}{1}{--}
    \subdividedEdge{$(gn0)+(0.5,0)$}{$(g6n0)+(0.5,0)$}{1}{to[bend angle=15, bend left] node[ivertex, pos=0.5]{}}
    \leafEdge{$(gn0)+(0.5,0)$}{u0}{}{above right=-0.05cm}{--}
    \leafEdge{$(gn0)+(-0.5,0)$}{v0}{7}{below left=-0.05cm}{--}
    \node[ivertex] at ($(gn0)+(-0.5,0)$) {};
    \node[ivertex] at ($(gn0)+(0.5,0)$) {};
    \node[above] at ($(gn0)+(-0.5,0)$) {$z_0$};
    \subdividedEdge{g1n0}{g1n1}{1}{--}
    \subdividedEdge{g1n1}{$(g9n0)+(0.5,0)$}{1}{.. controls (g1n1-|v13) and (g9n0-|v3) ..}
    \subdividedEdge{g1n0}{g4n0}{1}{to[bend angle=35, bend right]}
    \subdividedEdge{g1n1}{$(g5n0)+(-0.5,0)$}{1}{--}
    \leafEdge{g1n0}{u1}{3}{left=-0.05cm}{--}
    \node[ivertex] at (g1n0) {};
    \node[ivertex] at (g1n1) {};
    \node[above] at (g1n1) {$z_1$};
    \subdividedEdge{$(g6n0)+(0.5,0)$}{$(g7n0)+(-0.5,0)$}{1}{to[bend angle=15, bend left]}
    \subdividedEdge{$(g6n0)+(0.5,0)$}{$(g8n0)+(-0.5,0)$}{1}{to[bend angle=45, bend left]}
    \leafEdge{$(g6n0)+(0.5,0)$}{u2}{4}{above left=-0.05cm}{--}
    \leafEdge{$(g9n0)+(0.5,0)$}{v1}{7}{right=-0.05cm}{to[bend angle=25, bend left]}
    \leafEdge{$(g9n0)+(0.5,0)$}{u3}{5}{above right=-0.05cm}{--}
    \leafEdge{g4n0}{u4}{5}{right=-0.05cm}{--}
    \leafEdge{g4n0}{v2}{7}{above}{--}
    \leafEdge{g4n0}{v3}{7}{right=-0.05cm}{--}
    \leafEdge{g4n0}{v4}{7}{above}{--}
    \node[ivertex] at (g4n0) {};
    \node[below=0.1cm] at ($(g4n0)+(-0.2,0)$) {$z_4$};
    \subdividedEdge{$(g5n0)+(-0.5,0)$}{$(g5n0)+(0.5,0)$}{1}{--}
    \leafEdge{$(g5n0)+(0.5,0)$}{u5}{5}{above right=-0.05cm}{--}
    \leafEdge{$(g5n0)+(-0.5,0)$}{v5}{6}{above left=-0.05cm}{--}
    \leafEdge{$(g5n0)+(0.5,0)$}{v6}{6}{above right=-0.05cm}{--}
    \node[ivertex] at ($(g5n0)+(-0.5,0)$) {};
    \node[ivertex] at ($(g5n0)+(0.5,0)$) {};
    \node[below=0.1cm] at ($(g5n0)+(-0.5,0)$) {$z_5$};
    \subdividedEdge{$(g6n0)+(-0.5,0)$}{$(g6n0)+(0.5,0)$}{1}{--}
    \leafEdge{$(g6n0)+(-0.5,0)$}{u6}{5}{above left=-0.05cm}{--}
    \leafEdge{$(g6n0)+(-0.5,0)$}{v7}{6}{above left=-0.05cm}{--}
    \leafEdge{$(g6n0)+(0.5,0)$}{v8}{6}{above right=-0.05cm}{--}
    \node[ivertex] at ($(g6n0)+(-0.5,0)$) {};
    \node[ivertex] at ($(g6n0)+(0.5,0)$) {};
    \node[left=0.1cm] at ($(g6n0)+(0.85,-0.35)$) {$z_2, z_6$};
    \subdividedEdge{$(g7n0)+(-0.5,0)$}{$(g7n0)+(0.5,0)$}{1}{--}
    \leafEdge{$(g7n0)+(0.5,0)$}{u7}{5}{above right=-0.05cm}{--}
    \leafEdge{$(g7n0)+(-0.5,0)$}{v9}{6}{above left=-0.05cm}{--}
    \leafEdge{$(g7n0)+(0.5,0)$}{v10}{6}{above right=-0.05cm}{--}
    \node[ivertex] at ($(g7n0)+(-0.5,0)$) {};
    \node[ivertex] at ($(g7n0)+(0.5,0)$) {};
    \node[below=0.1cm] at ($(g7n0)+(-0.5,0)$) {$z_7$};
    \subdividedEdge{$(g8n0)+(-0.5,0)$}{$(g8n0)+(0.5,0)$}{1}{--}
    \leafEdge{$(g8n0)+(0.5,0)$}{u8}{5}{above right=-0.05cm}{--}
    \leafEdge{$(g8n0)+(-0.5,0)$}{v11}{6}{above left=-0.05cm}{--}
    \leafEdge{$(g8n0)+(0.5,0)$}{v12}{6}{above right=-0.05cm}{--}
    \node[ivertex] at ($(g8n0)+(-0.5,0)$) {};
    \node[ivertex] at ($(g8n0)+(0.5,0)$) {};
    \node[below=0.1cm] at ($(g8n0)+(-0.5,0)$) {$z_8$};
    \subdividedEdge{$(g9n0)+(-0.5,0)$}{$(g9n0)+(0.5,0)$}{1}{--}
    \leafEdge{$(g9n0)+(-0.5,0)$}{u9}{5}{above left=-0.05cm}{--}
    \leafEdge{$(g9n0)+(-0.5,0)$}{v13}{6}{above left=-0.05cm}{--}
    \leafEdge{$(g9n0)+(0.5,0)$}{v14}{6}{above right=-0.05cm}{--}
    \node[ivertex] at ($(g9n0)+(-0.5,0)$) {};
    \node[ivertex] at ($(g9n0)+(0.5,0)$) {};
    \node[left=0.1cm] at ($(g9n0)+(0.85,-0.35)$) {$z_3, z_9$};
    \foreach \u in {0,1,...,9}{
      \node[vertex] at (u\u) {};
      \node[above] at (u\u) {$u_{\u}$};
    }
    \foreach \v in {0,1,...,14}{\node[vertex] at (v\v) {};}
    \foreach \v in {2,3,...,14}{\node[below] at (v\v) {$v_{\v}$};}
    \node[below] at (v0) {$v_{0}$};
    \node[above] at (v1) {$v_{1}$};
    \node at (l0) {$T'$};
    \node at (l1) {$T'_1$};
    \node at (l2) {$T'_2$};
    \node at (l3) {$T'_3$};
    \node at (l4) {$T'_4$};
    \node at (l5) {$T'_5$};
    \node at (l6) {$T'_6$};
    \node at (l7) {$T'_7$};
    \node at (l8) {$T'_8$};
    \node at (l9) {$T'_9$};
  \end{tikzpicture}\end{center}
  \caption{A {\TPG} $G$ (top), the cotree $\cT$ of $G$ (2nd row), an optimal (odd) $11$-leaf root $T$ of $G$ (3rd row) as computed by Algorithm~\ref{alg:opt} (with input $p=1$), and an optimal even $12$-leaf root $T'$ of $G$ (bottom) as computed by Algorithm~\ref{alg:opt} with input $p = 0$.}\label{fig:algorithm}
\end{figure}
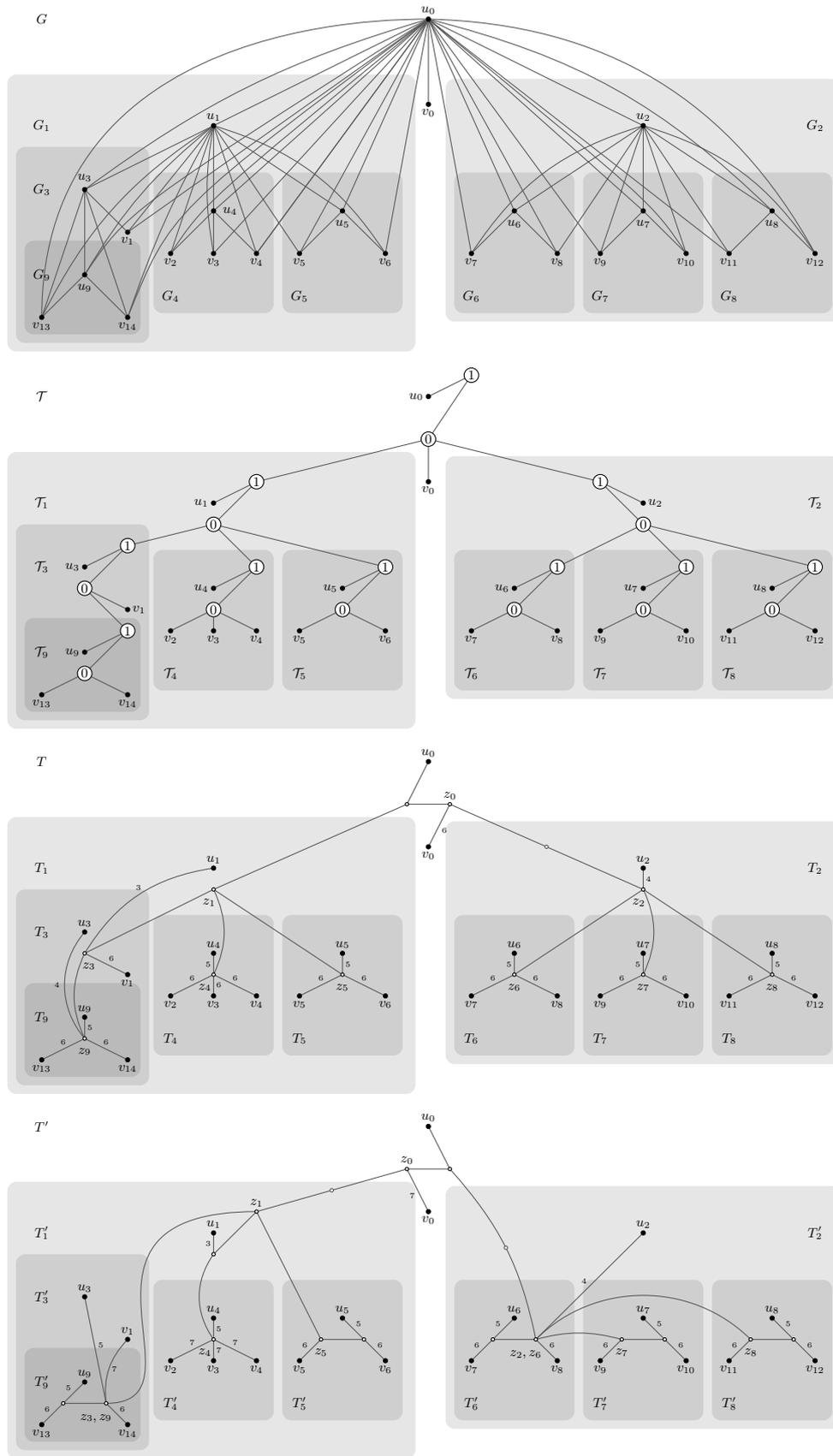

\smallskip\noindent
Hence, if the input graph $G$ is not a star then the approach is to firstly divide $G$ into smaller connected subgraphs $G_1, \dots, G_s$ (and isolated vertices $v_1, \dots, v_t$), secondly find corresponding $k$-leaf roots $T_1, \dots, T_s$ by recursion and the $\eta$-operation, and, last, conquer by merging them into a single leaf root of $G$ with the $\mu$-operation.
The divide-step is simple for disconnected $G$ and, otherwise, is carried out by removing the unique universal (cutvertex) of $G$.

The $\rho$-operation is sensitive to the given parity $p$ for using $\mu$ as a subroutine.
Observe that $p$ also decides how the resulting $k$ is determined.
There are four cases when $G$ is connected and not a star.
In the first one, when $s = 1$, the construction is the same for odd and even $p$ and consists of adding $u, v_1, \dots, v_t$ at the correct distance to the center of $T_1$ and computing $k$ from $k_1$.
Secondly, if $s > 1$ and $p = 1$, the $\mu$-operation has only one way of merging the recursively found leaf roots $T_1, \dots, T_s$ to minimize the diameter of the result $T$.
Then, $k$ widely depends on the two largest values of $k_1, \dots, k_s$.
But if $p = 0$, there is one situation that, on the one hand, allows $\mu$ to use a smaller diameter for $T$ by prioritizing the critical leaf root and, on the other hand, lets $\rho$ return a slightly better value for $k$.
This happens only when $s = 2$, or whenever the $k_c$-leaf root, with $k_c$ the third-largest value among $k_1, \dots, k_s$, properly fits into the diametral space of $T$ that is already required for the $k_a$-leaf root and the $k_b$-leaf root.

The third and bottom row of Figure~\ref{fig:algorithm} illustrate the results $(T, 11)$ of $\rho(G, 1)$ and $(T', 12)$ of $\rho(G, 0)$ on the example $G$.
By recursion, both are produced bottom-up, and it is difficult to follow their assembly at the deeper recursion levels.
The highest recursion level of $\rho(G,1)$, however, has received a $7$-leaf root with odd diameter for subgraph $G_1$ and a $5$-leaf root with even diameter for $G_2$ in Step~(\ref{const:case:connected}\ref{const:case:connected:recursion}).
In Step~(\ref{const:case:connected}\ref{const:case:connected:k}), the $\rho$-procedure determines $k = k_1 + k_2 - 1 = 11$.
The extension of the trees in Step~(\ref{const:case:connected}\ref{const:case:connected:widen}) produces the $11$-leaf roots $T_1$ and $T_2$ for $G_1$ and $G_2$, respectively, as shown in Figure~\ref{fig:algorithm}. 
Their following merging and the attachment of $u_0$ in Step~(\ref{const:case:connected}\ref{const:case:connected:assembly}) produces the shown $11$-leaf root $T$ of $G$.
Similarly, $\rho(G,0)$ receives an odd-diameter $8$-leaf root of $G_1$ and an even-diameter $6$-leaf root of $G_2$.
Since $s = 2$, the critical root can be treated in the special way and, thus, $\rho(G,0)$ determines $k' = k'_1 + k'_2 - 2 = 12$.
After the extension, we get the $12$-leaf roots $T'_1$ for $G_1$ and $T'_2$ for $G_2$ as in Figure~\ref{fig:algorithm}.
Their merging and the attachment of $u_0$ yields the $12$-leaf root $T'$ of $G$ as also illustrated there.

The following statement regards the correctness of our procedure.
\begin{theorem}\label{cor:correct_optimal}
  Let $G$ be a {\TPG} on $n$ vertices and without true twins and let $p \in \{0,1\}$.
  Then $(T,k) = \rho(G,p)$ provides a $k$-leaf root $T$ of $G$ that is optimal for parity $p$ (hence, $\odd(k) = p$) and with $k \leq n+1$.
  If $G$ is connected then
  \begin{description}[nosep]
    \item[\em (T1)] $\rad(T) = k-1$,
    \item[\em (T2)] $\dmin_T = 1+\odd(T)$, and
    \item[\em (T3)] $\diam(T')\ge \diam(T)+k'-k$ for all $k'$-leaf roots $T'$ of $G$ with $\odd(k')=p$.
  \end{description}
\end{theorem}
Note that, with respect to the optimality of the result, the theorem above makes a slightly stronger statement than our main result in Theorem~\ref{thm:main}.
In fact, the $\rho$-operation can find a $\kappa$-leaf root with minimum $\kappa$ for every {\TPG} $G$ simply by choosing the best from $(T,k) = \rho(G, 1)$ and $(T',k') = \rho(G, 0)$.
To prove Theorem~\ref{thm:main}, the next section shows how to implement the $\rho$-operation in linear time.

\section{Linear Time Leaf Root Construction for CCGs}
\label{sec:cotree}
  
The algorithm in this section is an implementation of the $\rho$-operation from Section~\ref{sec:tplr}.
Here, the recursive subdivision of the input {\TPG} $G$ is replaced with a post-order iteration of the cotree of $G$.
But before we go into the details, we analyze the used submodules and show that the operations $\eta$ and $\mu$ run efficiently.
\begin{lemma}\label{lem:subdividing:complexity}
  Let $T$ be a compressed tree with $n$ leaves and with explicitly given min-max center $z$, center $Z$, diameter $\diam(T)$, and leaf-distance $\dmin_T$.
  For all integers $\delta \geq 0$, the compressed tree $T' = \eta(T,\delta)$ with min-max center $z'$, center $Z'$, diameter $\diam(T')$, and leaf distance $\dmin_{T'}$ can be computed in $\cO(n)$ time.
\end{lemma}
\begin{lemma}\label{lem:disconnected:complexity}
  Let $s \geq 0$ be an integer and, for all $i \in \{1, \dots, s\}$, let $T_i$ be a given, compressed tree with explicitly given min-max center $z_i$, center $Z_i$, diameter $\diam(T_i)$, and leaf-distance $\dmin_{T_i}$.
  For all integers $k \geq 2$ and vertices $v_1, \dots, v_t$, the merged compressed tree $T' = \mu(k, T_1, \dots, T_s, v_1, \dots, v_t)$ with center $Z'$ and diameter $\diam(T')$ can be computed in $\cO(s + t)$ time.
\end{lemma}
The recursive definition of $\rho(G,p)$ is implemented as an iterative traversal of the cotree of $G$.
We observe that connected and disconnected graphs $G$ can easily be distinguished with the cotree of $G$;
the former have a root labelled by $\otimes$ and the latter by  $\oplus$.
Likewise, we detect small input graphs, stars, with the cotree by checking if the root is labelled with $\otimes$ and if the only child that is labelled with $\oplus$ has just leaf-children.

Recall that, for connected graphs $G$ of sufficient size, the $\rho$-operation divides $G$ at the unique universal vertex $u$, to recurse into the non-trivial connected components $G_1, \dots, G_s$ of $G-u$, and to conquer by merging the according leaf roots $T_1, \dots, T_s$ into a parity-optimal solution for $G$.
This divide-and-conquer procedure is translated into a traversal of the cotree $\cT$ as follows.
Since input consists of {\TPGs} without true twins, we rely on Proposition~\ref{lem:trivially_perfect:cotree}.
That means that nodes with the label $\otimes$, like the root $X$ of $\cT$, always have exactly one leaf-child, say $u$, and one child with label $\oplus$, say $Y$.
The leaf $u$ marks the unique universal vertex in $G$ and $Y$ has children $Z_1, \dots, Z_{s}$ with label $\otimes$ and leaf-children $v_1, \dots v_t$ that represent the non-trivial connected components $G_1, \dots, G_s$ and the isolated vertices $v_1, \dots, v_t$ of $G-u$.
The chosen post-order traversal of $\cT$ makes sure that, before processing $X$ (and $Y$), the nodes $Z_1, \dots, Z_{s}$ have been visited and finished.
Because we use a stack to pass interim results upwards, we always find leaf roots $T_1, \dots, T_s$ for $G_1, \dots, G_s$ on the stack (in reverse order), when we need to compute a leaf root $T$ for the subgraph that corresponds to $\cT_X$.

We present the details of our construction in Algorithm~\ref{alg:opt}:~\texttt{OptimalLeafRoot} and summarize our results in the following theorem.
\begin{theorem}\label{thm:main:odd_even}
  Given a chordal cograph $G$ on $n$ vertices and $m$ edges and $p \in \{0,1\}$, a (compressed) $\kappa$-leaf root of $G$ with minimum integer $\kappa$ of parity $p$ can be computed in $\cO(n + m)$ time.
\end{theorem}

\begin{algorithm}[h!]
  \DontPrintSemicolon 
  \KwIn{A {\TPG} $G=(V,E)$ without true twins and a parity $p \in \{0,1\}$}
  \KwOut{A pair $(T, k)$ of a $k$-leaf root $T$ of $G$ with smallest $p$-parity integer $k$.}
  
  \smallskip
  initialize empty stack $\cS$ and compute the cotree $\cT$ of $G$\;\label{algo:init_ready}
  \ForEach{node $X$ visited traversing $\cT$ in post-order}{\label{algo:loop_start}
    \If{$X$ is labelled with $\otimes$}{
      let $Y$ be the $\oplus$-child and $u$ the leaf-child of $X$\;
      let $s$ be the number of $\otimes$-children and $v_1, \dots, v_t$ the leaf-children of $Y$\;
      \uIf(\tcp*[f]{Case~\ref{const:base}, base case, input is a star}){$s = 0$}{
        build $T$ like Case~\ref{const:base}, Section~\ref{sec:tplr} for star on edges $uv_1, \dots, uv_t$\;\label{algo:start_base_case}
        push $(T, 4-p)$ onto $\cS$\;\label{algo:base_push}
      }
      \Else(\tcp*[f]{Case~\ref{const:case:connected}, input is a connected graph}){
        \lForEach{$i\in \{s, s-1, \dots, 1\}$}{pop $(T_i,k_i)$ from $\cS$}\label{algo:start_connected_case}
        \lIf{$s = 1$}{$k \gets k_1 + 2(1 - \odd(T_1))$}
        \Else{
          $k_a \gets \max\{k_1, \dots, k_s\}$\;
          $k_b \gets \max\{k_i \mid 1 \leq i \leq s, i \not= a\}$\;
          \lIf{$p = 1$}{$k \gets k_a + k_b - 1 - 2\cdot\odd(T_a)\cdot\odd(T_b)$}
          \Else{
            \uIf{$s > 2$ and $k_a > \max\{k_i \mid 1 \leq i \leq s, i \not= a, i \not= b\}$}{
              $k \gets k_a + k_b - 2\cdot\odd(T_a)\cdot\odd(T_b)$\;
            }
            \lElse{$k \gets k_a + k_b - 2\cdot(\odd(T_a)+\odd(T_b)-\odd(T_a)\cdot\odd(T_b))$}
          }
        }
        \lForEach{$i\in \{1, \dots, s\}$}{$T'_i \gets \eta(T_i, (k-k_i)/2)$}\label{algo:connected_widening}
        $T \gets \mu(k, T'_1, \dots, T'_s, v_1, \dots, v_t)$ with $ \gets$ a center of $T$\;\label{algo:connected_merging}
        attach $u$ as a leaf to a center of $T$\;\label{algo:connected_attach}
        push $(T,k)$ onto $\cS$\;\label{algo:connected_push}
      }
    }\label{algo:loop_end} 
    \If(\tcp*[f]{Case~\ref{const:case:disconnected}, disconnected input}){$X$ is $\oplus$-node without parent}{
      let $s$ be the number of $\otimes$-children and $v_1, \dots, v_t$ the leaf-children of $X$\;\label{algo:disconnected_start}
      \lForEach{$i\in \{1, \dots, s\}$}{pop $(T_i,k_i)$ from $\cS$}
      $k \gets \max\{k_1, \dots, k_s, p+2\}$\;
      \lForEach{$i\in \{s, s-1, \dots, 1\}$}{$T'_i \gets \eta(T_i, (k-k_i)/2)$}\label{algo:disconnected_widening}
      $T \gets \mu(k, T'_1, \dots, T'_s, v_1, \dots, v_t)$\;\label{algo:disconnected_merging}
      push $(T,k)$ onto $\cS$\;\label{algo:disconnected_push}
    }
    pop $(T,k)$ from $\cS$\;\label{algo:interim_result}
    \Return $(T,k)$\;\label{algo:result}
  }
  \caption{\texttt{OptimalLeafRoot}}
  \label{alg:opt}
\end{algorithm}

\section{Conclusion}
\label{sec:oddeven}

With Theorem~\ref{thm:main:odd_even}, we have shown that the {\OLR} problem is linear-time solvable for chordal cographs.
Our work also provides a linear-time solution for the $k$-leaf power recognition problem on chordal cographs.
Specifically, for a given {\TPG} $G$ and an integer $k$, it is sufficient to compute $(T, \kappa) = \rho(G, \odd(k))$ (in linear time with Algorithm~\ref{alg:opt}) to see by $\kappa \leq k$ if $G$ is a $k$-leaf power.

\noindent
We conclude the paper by exploring the differences in the construction of odd and even leaf-roots.
As we have seen, merging three or more even leaf roots sometimes requires a slightly stronger increase in $k$ than for odd leaf roots.
This can accumulate to an arbitrary big gap between $k$ and $k'$ of an optimal odd $k$-leaf root and an optimal even $k'$-leaf root of a given {\TPG}.
For example, consider the (infinite) series $F_1, F_2, F_3, \dots, F_i, \dots$ of {\TPGs} defined as follows.
Let $F_0$ be the path on three vertices and for all integers $i > 0$ define
\[F_i = t_i\ \otimes\ ((x_i\ \otimes\ (F_{i-1}\ \oplus\ u_i))\ \oplus\ (y_i\ \otimes\ (F_{i-1}\ \oplus\ v_i))\ \oplus\ (z_i\ \otimes\ (F_{i-1}\ \oplus\ w_i)))\]
with $g \in \{t_i, u_i, v_i, w_i, x_i, y_i, z_i\}$ denoting a graph with the single vertex $g$.
By Section~\ref{sec:basic}, $F_1, F_2, \dots$ are a family of {\TPGs} and, apparently, all these graphs are connected and without true twins.
\begin{lemma}\label{lem:family_roots}
  For all integers $i \geq 1$, the graph $F_i$ is a (odd) $k_i$-leaf power for $k_i = 2^{i+2}-1$ but not a $(k_i-2)$-leaf power and a (even) $k'_i$-leaf power for $k'_i = k_i + 2^{i} - 1$ but not a $(k'_i-2)$-leaf power.
\end{lemma}
This means that, although odd and even leaf root construction follows the same approach, there are $k$-leaf powers of odd $k$ among the chordal cographs that have optimal even $k'$-leaf roots with $k'$ roughly $1.25k$.

\section{Proofs}
\label{sec:proofs}

\newcommand{\repetition}[3]{
\noindent
\begin{minipage}[t]{\textwidth}
\textbf{#1~\ref{#2}}. \emph{#3}
\end{minipage}
}


\repetition{Proposition}{lem:trivially_perfect:cotree}{
  If $G$ is a {\TPG} without true twins and $\cT$ is the cotree of $G$ then every $\otimes$-node of $\cT$ has exactly two children, one leaf and one $\oplus$-node. 
}
\begin{proof}
  Consider any $\otimes$-node $X$ of $\cT$.
  Then $X$ is not a leaf and, thus, has at least two children none of which a $\otimes$-node.
  Since $G$ is free of true twins, there is at most one leaf among the children of $X$.
  Distinct leaves $x, y$ would be true twins since, by definition, both are adjacent to all $z \in V(G)$ with the least common ancestor in $\cal T$ labelled by $\otimes$.
  Because $G$ is free of cycles on four vertices, $X$ also has at most one $\oplus$-child.
  Two distinct $\oplus$-nodes $Y, Z$ would represent two induced subgraphs $G_Y$ and $G_Z$, where, by definition, $G_Y$ contained two not adjacent vertices $a,b$ and $G_Z$ two not adjacent vertices $c,d$.
  Also by definition, there would be edges $ac,ad,bc,bd$ in $E(G)$ and, thus, an induced cycle on four vertices in $G$.
  Hence, $X$ having exactly two children, one leaf and one $\oplus$-node, is the only remaining possibility.
\qed\end{proof}

\repetition{Proposition}{lem:wolk_plus}{
  Every connected {\TPG} $G$ without true twins has a unique universal vertex $u$ and $G-u$ is disconnected (that is, $u$ is a cutvertex).
}
\begin{proof}
  Because $G$ is connected, the cotree has a $\otimes$-root $X$.
  By Proposition~\ref{lem:trivially_perfect:cotree}, $X$ has a leaf child $u$ that, by definition, is universal in $G$ and a $\oplus$-child $Y$.
  Additional universal vertices would be true twins of $u$.
  By definition, $G-u = G_Y$ is disconnected and, thus, $u$ is a cutvertex.
\qed\end{proof}

\repetition{Lemma}{lem:technical}{
  If $T_1, \dots, T_s$ are $s \ge 2$ trees with $\diam(T_1) \ge \dots \ge \diam(T_s)$ then
  \begin{enumerate}[label=\textup{(\roman*)},leftmargin=*]
  \item
  $\rad(T_1)\ge \rad(T_2)\ge \cdots \ge \rad(T_s)$.
  \item
  For all $1\le i<j\le s$, if $\rad(T_i)=\rad(T_j)$ then $\odd(T_i)\le \odd(T_j)$.
  \item
  For all $1\le i<j \le s$, $\rad(T_i)-\odd(T_i)\ge \rad(T_j)-\odd(T_j)$. 
  \end{enumerate}
}
\begin{proof}
  Let $i<j$.
  \begin{enumerate}[label=\textup{(\roman*)},leftmargin=*]
  \item[\ref{lem:technical:i}]
  By definition, $2\rad(T_i)-\odd(T_i) = \diam(T_i) \geq \diam(T_j) = 2\rad(T_j)-\odd(T_j)$.
  Clearly, if $\odd(T_i) = \odd(T_j)$ then we directly get $\rad(T_i) \geq \rad(T_j)$.
  Otherwise, if only $\odd(T_i) = 1$, we argue that
  \[2\rad(T_i) > 2\rad(T_i)-1 = \diam(T_i) \geq \diam(T_j) = 2\rad(T_j).\]
  Finally, if just $\odd(T_j) = 1$ then the diameter of $T_i$ is even and the diameter of $T_j$ is odd, thus, not equal.
  This immediately tells us that
  \[2\rad(T_i) = \diam(T_i) > \diam(T_j) = 2\rad(T_j)-1\]
  and, hence, $2\rad(T_i) \geq 2\rad(T_j)$ and we are done.
  \item[\ref{lem:technical:ii}]
  If $\rad(T_i)=\rad(T_j)$, then $2\rad(T_i)-\odd(T_i)=\diam(T_i)\ge\diam(T_j)=2\rad(T_j)-\odd(T_j)=2\rad(T_i)-\odd(T_j)$.
  Hence, $\odd(T_i)\le \odd(T_j)$.
  \item[\ref{lem:technical:iii}]
  By \ref{lem:technical:i}, $\rad(T_i)\ge \rad(T_j)$.  
  If $\rad(T_i)>\rad(T_j)$, then $\rad(T_i)-\odd(T_i)> \rad(T_j)-\odd(T_i)\ge \rad(T_j)-1$.
  Hence, $\rad(T_i)-\odd(T_i)\ge \rad(T_j)\ge \rad(T_j)-\odd(T_j)$.  
  If $\rad(T_i)=\rad(T_j)$, then, by \ref{lem:technical:ii}, $\odd(T_i)\le \odd(T_j)$.
  Hence, $\rad(T_i)-\odd(T_i)=\rad(T_j)-\odd(T_i)\ge \rad(T_j)-\odd(T_j)$.
  \end{enumerate}
\qed\end{proof}

\repetition{Lemma}{lem:dmin}{
  Every $k$-leaf root $T$ of a connected graph satisfies $\dmin_T \le \tfrac{k}{2}$. 
}
\begin{proof}
  Let $T$ be a $k$-leaf root of a connected graph $G$, and let $z$ be a min-max center vertex of $T$.
  Let $T'$ be any subtree of $T-z$.
  Then, every pair of leaf $v$ in $T'$ and leaf $w$ outside $T'$ fulfills
  \[\dist_T(u,v) = \dist_T(u,z) + \dist_T(z,v) \ge 2\!\cdot\! \dmin_T.\]
  Thus, if $\dmin_T$ was larger than $\tfrac{k}{2}$, then $\dist_T(u,v) > k$ and $G$ would be disconnected.
\qed\end{proof}

\repetition{Lemma}{lem:A}{
  If $G$ is a non-complete graph with a universal vertex $u$ and $T$ is a $k$-leaf root of $G$ then $\dist_T(u,z) \le k-\rad(T)+\odd(T)$ for all center vertices $z$ of $T$.
  If $z_1 \not= z_2$ are the center vertices of $T$, then $\dist_T(u,z_1) \leq k-\rad(T)$ or $\dist_T(u,z_2) \leq k-\rad(T)$.
}
\begin{proof}
  Consider a diametral path $P$ in $T$, and let $x, y$ be the end vertices of $P$.
  Since $G$ is not complete, $x$ and $y$ are not adjacent in $G$.
  Hence, $u$, $x$ and $y$ are different leaves in $T$. 
  Then the three paths in $T$, the $x,u$-path, the $y,u$-path, and $P$ intersect at a unique (non-leaf) vertex on $P$, say $c$. 
  Let $\ell_1$ be the length of the $x,c$-path and $\ell_2$ be the length of the $y,c$-path in $T$. 
  Since $u$ is universal in $G$, we have
  \[\dist_T(u,c)+\ell_1\le k\quad \text{and}\quad \dist_T(u,c)+\ell_2\le k.\]
  Let $z_1, z_2$ be the center vertices of $T$; possibly with $z_1 = z_2$, if $\diam(T)$ is even. 
  Recall that, since $P$ is a diametral path, the centers of $T$ and $P$ coincide, and that $z_1z_2$ is an edge on $P$ whenever $z_1 \not= z_2$.
  Let, without loss of generality, $z_1$ be on the $x,z_2$-path (and so $z_2$ is on the $y,z_1$-path). 
  Then
  \begin{align*}
    \dist_T(z_1,y) &= \rad(T) = \dist_T(z_2,x) \text{ and}\\
    \dist_T(z_1,x) &= \rad(T)-\odd(T) = \dist_T(z_2,y).
  \end{align*} 
  Assume $c$ on the $x,z_1$-path; possibly with $c=z_1$.
  Then
  \begin{align*}
  \dist_T(u,z_1) & = \dist_T(u,c)+\big(\ell_2-\dist_T(z_1,y)\big)\\
                & = \big(\dist_T(u,c)+\ell_2\big)-\rad(T)\\
                & \le k-\rad(T) \quad\text{and}\\
  \dist_T(u,z_2) & = \dist_T(u,c)+\big(\ell_2-\dist_T(z_2,y)\big)\\ 
                & = \big(\dist_T(u,c)+\ell_2\big)-\big(\rad(T)-\odd(T)\big)\\
                & \le k-\rad(T)+\odd(T). 
  \end{align*} 
  In the case, where $c$ is on the $y,z_2$-path, the same arguments yield
  $\dist_T(u,z_2) \le k-\rad(T)$ and $\dist_T(u,z_1) \le k-\rad(T)+\odd(T)$.
\qed\end{proof}

\repetition{Corollary}{cor:A}{
  If $G$ is a graph with a universal vertex and $T$ is a $k$-leaf root of $G$ then $\rad(T)\le k-1$ and, in particular, $\diam(T)\le 2k-2$.
}
\begin{proof}
  For complete graphs, the statement is obvious.
  So, let $G$ be a non-complete graph with universal vertex $u$.
  By Lemma~\ref{lem:A},
  $1 \le \dist_T(u,z) \le k-\rad(T)$
  for any center vertex $z$ of $T$.
  Hence, $\rad(T) \le k-1$.
  Since $\diam(T)\le 2\!\cdot\!\rad(T)$, we have $\diam(T)\le 2k-2$.
\qed\end{proof}


\repetition{Lemma}{lem:subdividing}{
  If $T$ is a $k$-leaf root of a graph $G$ and $\delta \ge 0$ an integer then $T' = \eta(T, \delta)$ is a $(k+2\delta)$-leaf root of $G$ with same center, same min-max center vertices, and
  $\diam(T') = \diam(T)+2\delta$, $\rad(T') = \rad(T)+\delta$, $\dmin_{T'}=\dmin_{T}+\delta$.
}
\begin{proof}
  Note that $T$ and $T'$ have the same set of leaves.
  For every two leaves $u$ and $v$, we have $\dist_{T'}(u,v)=\dist_{T}(u,v)+2\delta$.
  Hence, $\dist_{T'}(u,v)\le k+2\delta$ if and only if $\dist_T(u,v)\le k$. 
  Thus, $T'$ is a $k+2\delta$-leaf root of $G$.
  The other statements are obvious from the respective definitions.
\qed\end{proof}
  
\repetition{Lemma}{lem:disconnected}{
  Let $G$ be a graph and $k \geq 2$ an integer.
  If $G$ is disconnected with $s \ge 0$ non-trivial connected components $G_1, \dots, G_s$ and $t \ge 0$ isolated vertices $v_1, \dots, v_t$ such that $s+t \ge 2$ and if $T_1, \dots, T_s$ are $k$-leaf roots for $G_1, \dots, G_s$ then $T = \mu(k, T_1, \dots, T_s, v_1, \dots, v_t)$ is a $k$-leaf root of $G$.
}
\begin{proof}
  From the construction, we immediately see for all distinct vertices $x$ and $y$ from the same component $G_i$ that $\dist_T(x,y) = \dist_{T_i}(x,y)$.
  Hence, $\dist_T(x, y)\le k$, if and only if $\dist_{T_i}(x,y) \le k$, if and only if $xy$ is an edge in $G_i$.

  To prove that $T$ is a $k$-leaf root for $G$, it is sufficient to show for all vertices $x$ and $y$ stemming from different components of $G$ that $\dist_T(x,y) > k$.
  In case of $s > 0$, we assume, without loss of generality, that the critical index is $m = 1$.
  For the remainder of the proof, we keep in mind that, by the definition of a min-max center, $\dist_{T_i}(x,z_i)-\dmin_{T_i} \geq 0$ for all $T_i, 1 \leq i \leq s$ with min-max center $z_i$ and all leaves $x$ of $T_i$.

  For a start, we consider a vertex $x$ in $G_1$ and $v_j$ for any $j \in \{1, \dots, t\}$ and see
  \[\dist_T(x,v_j) = \dist_{T_1}(x, z_1) + \big(\tfrac{k+\odd(k)}{2} - \dmin_{T_1}\big) + \big(\tfrac{k-\odd(k)}{2} + 1\big) \ge k+1.\] 
  Similarly, for any vertex $x$ in $G_i$, $1 < i \le s$, and $v_j$, we get
  \[\dist_T(x,v_j) = \dist_{T_i}(x, z_i) + \big(\tfrac{k-\odd(k)}{2}+1-\dmin_{T_i}\big) + \big(\tfrac{k-\odd(k)}{2} + 1\big) \ge k+1.\]
  For any vertex $x$ in $G_1$ and $y$ in $G_j$ with $1 < j \leq s$, it is
  \begin{align*}
  \dist_T(x,y) &= \dist_{T_1}(x, z_1) + \big(\tfrac{k+\odd(k)}{2} - \dmin_{T_1}\big) + \big(\tfrac{k-\odd(k)}{2}+1-\dmin_{T_j}\big) + \dist_{T_j}(z_j, y)\\
                &= (k+1) + \big(\dist_{T_1}(x,z_1)-\dmin_{T_1}\big) + \big(\dist_{T_j}(z_j, y)-\dmin_{T_j}\big) \ge k+1.
  \end{align*}
  Similarly, for any vertex $x$ in $G_i$ and $y$ in $G_j$ with $1 < i<j \le s$, we get 
  \begin{align*}
  \dist_T(x,y) &= \dist_{T_i}(x, z_i) + \big(\tfrac{k-\odd(k)}{2}+1-\dmin_{T_i}\big) + \big(\tfrac{k-\odd(k)}{2}+1-\dmin_{T_j}\big) + \dist_{T_j}(z_j, y)\\
                &= (k+2-\odd(k)) + \big(\dist_{T_i}(x,z_i)-\dmin_{T_i}\big) + \big(\dist_{T_j}(z_j, y)-\dmin_{T_j}\big) \ge k+1. 
  \end{align*}
  Finally, the distance between $v_i$ and $v_j$ with $1 \le i<j \le t$ is
  \[\dist_T(v_i,v_j) = 2\cdot\big(\tfrac{k-\odd(k)}{2}+1\big) \ge k+1.\]
  It follows for all $x,y \in V(G)$ that $xy \in E(G)$ if and only if $\dist_T(x,y)\le k$.
  That is, $T$ is a $k$-leaf root of $G$.
\qed\end{proof}

\repetition{Theorem}{cor:correct_optimal}{
  Let $G$ be a {\TPG} on $n$ vertices and without true twins and let $p \in \{0,1\}$.
  Then $(T,k) = \rho(G,p)$ provides a $k$-leaf root $T$ of $G$ that is optimal for parity $p$ (hence, $\odd(k) = p$) and with $k \leq n+1$.
  If $G$ is connected then
  \begin{description}
    \item[\em (T1)] $\rad(T) = k-1$,
    \item[\em (T2)] $\dmin_T = 1+\odd(T)$, and
    \item[\em (T3)] $\diam(T')\ge \diam(T)+k'-k$ for all $k'$-leaf roots $T'$ of $G$ with $\odd(k')=p$.
  \end{description}
}

\begin{proof}
  The proof of Theorem~\ref{cor:correct_optimal} stands for the majority of the work behind this paper and, without subdivision, it is fairly long.
  So, to make reading easier, we break it down into the proofs of several propositions.
  In principle, the proof will be by complete induction on the vertex number and, below, we begin with the first proposition that serves the base case.
  
  \begin{proposition}\label{obs:star}
    Let $G$ be a star on $n \geq 3$ vertices and let $p \in \{0,1\}$ be a given parity.
    Then $(T,k) = \rho(G,p)$ provides a $k$-leaf root $T$ of $G$ that is optimal for parity $p$ (hence, $\odd(k) = p$) with $k = 4 - p \leq n+1$ and satisfying \textup{(T1)}, \textup{(T2)}, and \textup{(T3)}.
  \end{proposition}
  For a proof of Proposition~\ref{obs:star}, we start by observing that $G$ does not have true twins and, thus, $\rho(G,p)$ is well-defined.
  Moreover, since $G$ is not complete, every $3$-leaf root of $G$ is an optimal (odd) leaf root and every $4$-leaf root is an optimal even leaf root of $G$.
  By the same reason, any $k'$-leaf root $T'$ of $G$ satisfies $\diam(T')\ge k'+1$.
  
  We first consider the odd case, where $k = 3 = 4 - \odd(k)$, and show that the tree $T$ returned in Step~(\ref{const:base}\ref{const:base:odd}) is a $3$-leaf root of $G$ satisfying the claimed conditions.
  In fact, $T$ is just $\eta(G,1)$ with an additional leaf attached for $u$.
  This is obviously a $3$-leaf root of $G$ with $\rad(T)=2=k-1<n$, $\diam(T)=4$ and, thus, $\odd(T)=0$ and $\dmin_T = 1 = 1+\odd(T)$.
  Moreover, for any $k'$-leaf root $T'$ of $G$, we see that
  \[\diam(T') \ge k'+1 = 4+ k'-3 = \diam(T) + k'-k,\]
  which settles the odd case. 

  In the even case, we have $k=4$.
  If $G$ has exactly two leaves $v_1, v_2$ attached to the central vertex $u$ then Step~(\ref{const:base}\ref{const:base:even2}) returns the tree $T$ where the leaves $u, v_1$ and $v_2$ are attached to a vertex $v$ with paths of lengths one, two and three, respectively.
  Obviously, $T$ is a $4$-leaf root of $G$ with $\rad(T)=3=k-1=n$, $\diam(T)=5$ and, thus, $\odd(T)=1$ and $\dmin_T= 2 = 1+\odd(T)$.
  Every even $k'$-leaf root $T'$ of $G$ fulfills
  \[\diam(T') \ge k'+1 = 5+ k'-4 = \diam(T)+k'-k,\]
  which settles this case, too.
  
  Finally, let $p = 0$ and $t \geq 3$, where Step~(\ref{const:base}\ref{const:base:evenN}) constructs $T$ from $\eta(G, 2)$ by attaching a new leaf to $u$.
  Similar to the odd case, $T$ is obviously a $4$-leaf root of $G$ with $\rad(T)=3=k-1<n$, $\diam(T)=6$ and, thus, $\odd(T)=0$ and $\dmin_T=1 =1+\odd(T)$.
  To see \textup{(T3)}, suppose to the contrary that $\diam(T') < \diam(T) + k'-k = 6+k'-4 = k'+2$ for some even $k'$-leaf root $T'$ of $G$.
  Consider the vertex $w$ on the $v_1, v_2$-path in $T'$ that is closest to $v_3$ and let $d_i = \dist_{T'}(v_i,w)$ for all $i \in \{1, 2, 3\}$.
  Because of the diameter bound and since $v_1, v_2$ and $v_3$ are pairwise non-adjacent in $G$, we have $k'+1 \leq d_i+d_j < k'+2$ for all $1 \leq i<j \leq 3$.
  This means $d_1+d_2 = d_1+d_3 = d_2+d_3 = k'+1$.
  Hence, $2(d_1+d_2+d_3) = 3(k'+1)$, which implies that $3(k'+1)$ is even.
  This is a contradiction to the fact that $k'$ is even and, thus, the last case is settled.
  Thus, Proposition~\ref{obs:star} has been shown.

  \medskip
  With the following proposition, we only care about non-trivial connected input graphs and, for the moment, refrain from showing the optimality of $(T,k)$ or even property \textup{(T3)}.
  \begin{proposition}\label{lem:correct}
    Let $G$ be a connected {\TPG} on $n$ vertices and without true twins and let $p \in \{0,1\}$ be a given parity.
    Then $(T,k) = \rho(G,p)$ provides a $k$-leaf root $T$ of $G$ with $\odd(k) = p$ and $k \leq n+1$ and satisfying \textup{(T1)} and \textup{(T2)}.
  \end{proposition}
  The proof of Proposition~\ref{lem:correct} works by induction on the vertex number of $G$. 
  The smallest connected {\TPG} without true twins is the star with two leaves (since we ignore graphs with just one vertex). 
  In Proposition~\ref{obs:star}, this base case has already been settled.
  Moreover, Proposition~\ref{obs:star} allows to assume that $G$ is not a larger star in the following induction step.
    
  \medskip
  Next, let $G$ have more than three vertices.
  Because $G$ is connected and not a star, we are in Case~\ref{const:case:connected} of the procedure for $\rho$.
  Again because $G$ is connected and without true twins, it has a unique universal cutvertex by Proposition~\ref{lem:wolk_plus}.
  Hence, $H = G-u$ is disconnected.
  Since $G$ is not a star, $H$ has $s \ge 1$ non-trivial connected components $G_1, \ldots, G_s$ and $t\ge 0$ isolated vertices $v_1,\ldots,v_t$ such that $s+t \ge 2$.
  
  Note that, as for $G$, every $G_i$ is a connected {\TPG} without true twins but with fewer vertices than $G$.
  This means that the induction hypothesis holds and, thus, Step~(\ref{const:case:connected}\ref{const:case:connected:recursion}) provides a $k_i$-leaf root $T_i$ of parity $p$ satisfying \textup{(T1)} and \textup{(T2)} for every $G_i$, $1\le i\le s$.
  We observe that $k_i \geq 3$ for all $i \in \{1, \dots, s\}$ as none of $G_1, \dots, G_s$ is a complete graph (since they are neither isolated vertices nor contain true twins).
  
  At this point, we assume, without loss of generality, that $\diam(T_1) \geq \diam(T_2) \geq \dots \geq \diam(T_s)$.
  In addition to a simplified argumentation below, we get that
  \[k_1 \geq k_2 \geq \dots \geq k_s \text{ and also } k_1 - \odd(T_1) \geq k_2 - \odd(T_2) \geq \dots \geq k_s - \odd(T_s).\]
  This is firstly because of Lemma~\ref{lem:technical}, which implies
  \begin{align*}
  &\rad(T_1) \geq \rad(T_2) \geq \dots \geq \rad(T_s) \text{ and }\\
  &\rad(T_1) - \odd(T_1) \geq \rad(T_2)  - \odd(T_2) \geq \dots \geq \rad(T_s) - \odd(T_a),
  \end{align*}
  and secondly because of the induction hypothesis $k_i = \rad(T_i) + 1$, $1 \leq i \leq s$. 
  Hence, if $s > 1$ then Step~(\ref{const:case:connected}\ref{const:case:connected:k}) selects $k_a = k_1$, $k_b = k_2$ and, if it exists, $k_c = k_3$.
  
  Next, let $k$ be the number computed in Step~(\ref{const:case:connected}\ref{const:case:connected:k}) and note that $\odd(k) = p$ and $k > k_i, 1 \leq i \leq s$.
  By Lemma~\ref{lem:subdividing}, Step~(\ref{const:case:connected}\ref{const:case:connected:widen}) produces a $k$-leaf root $T'_i = \eta\left(T_i, \tfrac{k-k_i}{2}\right)$ of $G_i$ for all $i \in \{1, \dots, s\}$ such that
  \begin{align*}
    \diam(T'_i) &= \diam(T_i)+k-k_i,\\
    \rad(T'_i) &= \rad(T_i)+\tfrac{k-k_i}{2} = \tfrac{k+k_i}{2}-1,\\
    k-1 &\leq |V(G_i)| + (k-k_i), \text{ and}\\
    \dmin_{T'_i} &=\dmin_{T_i}+\tfrac{k-k_i}{2} = 1+\odd(T_i)+\tfrac{k-k_i}{2}.
  \end{align*}
  Then, Lemma~\ref{lem:disconnected} tells us that Step~(\ref{const:case:connected}\ref{const:case:connected:assembly}) assembles a $k$-leaf root $T'$ of $H$.
  
  Observe that, under the given assumption that $\diam(T_1) \geq \diam(T_2) \geq \dots \geq \diam(T_s)$, the $\mu$-operation chooses the critical index $m = 1$.
  This happens because $\mu$ uses $m$ to select the first tree in $T'_1, \dots, T'_s$ that minimizes $\dmin_{T'_m} = 1+\odd(T_m)+\tfrac{k-k_m}{2}$.
  In fact, if there was $i \in \{2, \dots, s\}$ with $\odd(T_i)-\tfrac{k_i}{2} < \odd(T_1)-\tfrac{k_1}{2}$ then we had $k_1 - k_i < 2(\odd(T_1)-\odd(T_i))$.
  Since $k_1 \geq k_i$ and both are of the same parity and $\odd(T_1), \odd(T_i) \in \{0,1\}$, we could conclude that $k_1 = k_i$ and $\odd(T_1) = 1$ and $\odd(T_i) = 0$.
  But that would imply the contradiction
  \begin{align*}
    \diam(T_1) &= 2\rad(T_1) - \odd(T_1) = 2(k_1-1) - 1 < 2(k_1-1) - 0\\
                &= 2(k_i-1) - \odd(T_i) = 2\rad(T_i) - \odd(T_i) = \diam(T_i).
  \end{align*}
  
  In the following, we show that $\rad(T') = k-1 \leq n$ and, for that purpose, we firstly determine the diameter of $T'$.
  A diametral path in $T'$ connects two leaves $x$ and $y$, thus, vertices of $H$.
  To find the length of a diametral path, thus, the diameter of $T'$, we subsequently analyze all possible origins of $x, y$ in $H$.

  For leaves $x$ and $y$ of the same $G_i, 1\leq i\leq s$, the longest possible connecting path has length
  \[
  \diam(T'_i) = 2\rad(T'_i) - \odd(T'_i) = 2\big(\tfrac{k+k_i}{2}-1\big) - \odd(T_i) = k+k_i- \odd(T_i)-2.
  \]
  The longest path in $T'$ connecting a leaf $x$ in $G_1$ (recall the critical index $m = 1$) and $y = v_j, 1 \leq j \leq t$ has length  
  \begin{align*}
  \ell_1 &= \big(\rad(T'_1) + \tfrac{k+\odd(k)}{2} - \dmin_{T'_1}\big) + \tfrac{k-\odd(k)}{2} + 1\\
        &= k + \rad(T'_1) - \dmin_{T'_1} + 1\\
        &= k + \big(\tfrac{k+k_1}{2} - 1\big) - \big(1+\odd(T_1)+\tfrac{k-k_1}{2}\big) + 1\\
        &= k + k_1 - \odd(T_1) - 1.
  \end{align*}
  Similarly, a longest path of $T'$ connecting a leaf $x$ in $G_i$ with $i \in \{2, \dots, s\}$ and $y = v_j, 1 \leq j \leq t$ has length
  \begin{align*}
  \ell_i &= \big(\rad(T'_i) + \tfrac{k-\odd(k)}{2}+1-\dmin_{T'_i}\big) + \tfrac{k-\odd(k)}{2} + 1\\
      &= k - \odd(k) + \rad(T'_i) - \dmin_{T'_i} + 1\\
      &= k - \odd(k) + \big(\tfrac{k+k_i}{2}-1) - \big(1+\odd(T_i)+\tfrac{k-k_i}{2}\big) + 1\\
      &= k - \odd(k) + k_i - \odd(T_i) - 1.
  \end{align*}
  Next, let $x$ be a leaf in $G_1$ (again, recall that the critical index $m$ is $1$) and $y$ a leaf in $G_i$ with $i \in \{2, \dots, s\}$ that are most distant from each other in $T'$. 
  The path connecting the leaves $x$ and $y$ has length
  \begin{align*}
  \ell_{1,i} &= \big(\rad(T'_1) + \tfrac{k-\odd(k)}{2}+1-\dmin_{T'_1}\big) + \big(\tfrac{k+\odd(k)}{2} - \dmin_{T'_i} + \rad(T'_i)\big)\\
      &= k + (\rad(T'_1)-\dmin_{T'_1}) + (\rad(T'_i)-\dmin_{T'_i}) + 1\\
      &= k + \tfrac{k+k_1}{2}-1 - \big(1+\odd(T_1)+\tfrac{k-k_1}{2}\big) + \tfrac{k+k_i}{2}-1 - \big(1+\odd(T_i)+\tfrac{k-k_i}{2}\big) + 1\\
      &= k + k_1 - \odd(T_1) + k_i - \odd(T_i) - 3.
  \end{align*}
  The leaves $x$ in $G_i$ and $y$ in $G_j$ with $i,j \in \{2, \dots, s\}$ and $i < j$ that are farthest from each other in $T'$ have distance
  \begin{align*}
  \ell_{i,j} &= \big(\rad(T'_i) + \tfrac{k-\odd(k)}{2}+1-\dmin_{T'_i}\big) + \big(\tfrac{k-\odd(k)}{2}+1-\dmin_{T'_j} + \rad(T'_j)\big)\\
      &= k-\odd(k) + (\rad(T'_i)-\dmin_{T'_i}) + (\rad(T'_j)-\dmin_{T'_j}) + 2\\
      &= k-\odd(k) + \tfrac{k+k_i}{2}-1 - \big(1+\odd(T_i)+\tfrac{k-k_i}{2}\big) + \tfrac{k+k_j}{2}-1 - \big(1+\odd(T_j)+\tfrac{k-k_i}{2}\big) + 2\\
      &= k-\odd(k) + k_i - \odd(T_i) + k_j - \odd(T_j) - 2.
  \end{align*}
  Finally, any $x = v_i$ and $y = v_j$, $1\le i<j \le t$, have distance
  \[\dist_{T'}(v_i, v_j) = 2\big(\tfrac{k-\odd(k)}{2} + 1\big) = k-\odd(k)+2.\]

  Before we can determine the diameter of $T'$, we need to sort out the longest paths of $T'$ from those given above.
  By premise, $\diam(T'_1) \geq \diam(T'_2) \geq \dots \geq \diam(T'_s)$.
  Next, if $t > 0$ then we obviously have that $\ell_1 > \diam(T'_1)$ and, by $k_1 \geq k_2 \geq \dots \geq k_s$, that $\ell_1 \geq \ell_i$ for all $i \in \{2, \dots, s\}$.
  Moreover, in case of $t \geq 2$, we also find that
  \[\ell_1 - \dist_{T'}(v_i, v_j) = (k + k_1 - \odd(T_1) - 1) - (k-\odd(k)+2) = k_1 + \odd(k) - \odd(T_1) - 3\]
  for all $1 \leq i < j \leq t$.
  If $k_1 = 3$ then, because $\odd(k) = \odd(k_1) = 1$, we get $\ell_1 - \dist_{T'}(v_i, v_j) = 1 - \odd(T_1) \geq 0$, which means $\ell_1 \geq \dist_{T'}(v_i, v_j)$.
  This holds even more so if $k_1 > 3$.

  If $s \geq 2$ then 
  \begin{eqnarray*}
    \ell_{1,2} - \ell_{1,i}  &=& (k + k_1 - \odd(T_1) + k_2 - \odd(T_2) - 3) - (k + k_1 - \odd(T_1) + k_i - \odd(T_i) - 3)\\
                              &=& ((k_2-1) - \odd(T_2)) - ((k_i-1) - \odd(T_i))\\
                              &=& (\rad(T_2)-\odd(T_2)) - (\rad(T_i)-\odd(T_i))\\
                              &\geq& 0  \text{ for all } i \geq 2 \text{ by Lemma~\ref{lem:technical}},\\  
    \ell_{1,2} - \diam(T'_1) &=& (k + k_1 - \odd(T_1) + k_2 - \odd(T_2) - 3) - (k + k_1 - \odd(T_1) - 2)\\
                              &=& k_2 - \odd(T_2) - 1 > 0 \text{ since } k_2 \geq 3, \text{ and}\\
      \text{if } t > 0 \text{ then } \ell_{1,2} - \ell_1 &=& (k + k_1 - \odd(T_1) + k_2 - \odd(T_2) - 3) - (k + k_1 - \odd(T_1) - 1)\\
                              &=& k_2 - \odd(T_2) - 2 \geq 0 \text{ because } k_2 \geq 3.
  \end{eqnarray*}
  Moreover, if $s > 2$ then
  \begin{eqnarray*}
    \ell_{1,2} - \ell_{2,3} &=& (k + k_1 - \odd(T_1) + k_2 - \odd(T_2) - 3) - (k-\odd(k) + k_2 - \odd(T_2) + k_3 - \odd(T_3) - 2)\\
                            &=& ((k_1-1)-\odd(T_1)) - ((k_3-1)-\odd(T_3)) + \odd(k) - 1\\
                            &=& (\rad(T_1)-\odd(T_1)) - (\rad(T_3)-\odd(T_3)) + \odd(k) - 1\\
                            &\geq& 0 \text{ if } \odd(k) = 1 \text{ by Lemma~\ref{lem:technical} and}\\
    \ell_{2,3} - \ell_{i,j} &=& (k-\odd(k) + k_2 - \odd(T_2) + k_3 - \odd(T_3) - 2) - (k-\odd(k) + k_i - \odd(T_i) + k_j - \odd(T_j) - 2)\\
                            &=& \big(((k_2-1)-\odd(T_2)) - ((k_i-1)-\odd(T_i))\big) + \big(((k_3-1)-\odd(T_3)) - ((k_j-1)-\odd(T_j))\big)\\
                            &=& \big((\rad(T_2)-\odd(T_2)) - (\rad(T_i)-\odd(T_i))\big) + \big((\rad(T_3)-\odd(T_3)) - (\rad(T_j)-\odd(T_j))\big)\\
                            &\geq& 0 \text{ if } 1 < i < j \leq s \text{ by Lemma~\ref{lem:technical}.}\\
  \end{eqnarray*}
  It remains to compare $\ell_{1,2}$ and $\ell_{2,3}$ in the even case where $\odd(k) = 0$:
  \begin{align*}
    &\ell_{1,2} - \ell_{2,3} = ((k_1-1)-\odd(T_1)) - ((k_3-1)-\odd(T_3)) + \odd(k) - 1\\
                            &= (k_1-\odd(T_1)) - (k_3- \odd(T_3)) - 1 \begin{cases}
                              \geq 0, &\text{if } k_1-\odd(T_1) > k_3-\odd(T_3),\\
                              < 0, &\text{otherwise}
                            \end{cases}
  \end{align*}
  
  We are now ready to estimate the diameter of $T'$ and, thus, the radius.
  If $s=1$ (and therefore $t > 0$) then
  \[\diam(T') = \max\{\diam(T'_1), \ell_1, \dots, \ell_t, \dist_{T'}(v_i, v_j) \mid 1\leq i < j \leq t\}.\]
  From the observations above, $\diam(T') = \ell_1$.
  In Step~(\ref{const:case:connected}\ref{const:case:connected:k}), we see that $k = k_1 + 2(1 - \odd(T_1))$ and so
  \[\ell_1 = k + k_1 - \odd(T_1) - 1 = k + (k - 2 + 2\odd(T_1)) - \odd(T_1) - 1 = 2k-3 + \odd(T_1).\]
  Obviously, $\odd(T') = 1-\odd(T_1)$, which means that $\rad(T') = k-1$.
  Note also that, since $n \geq |V(G_1)| + 2$ (because $\{u, v_1\} \subseteq V(G) \setminus V(G_1)$)
  \[k-1 \leq |V(G_1)| + (k-k_1) \leq (n - 2) + (2 - 2\odd(T_1)) \leq n.\]
  
  If $s \geq 2$ then let 
  \begin{align*}
    m_1 &= \max\{\diam(T'_1), \dots, \diam(T'_s)\},\\
    m_2 &= \max\{\ell_{i,j} \mid 1\leq i < j \leq s\},\\
    m_3 &= \max\{\ell_1, \dots, \ell_t\},\\
    m_4 &= \max\{\dist_{T'}(v_{i}, v_{j}) \mid 1 \leq i < j \leq t\}.
  \end{align*}
  Then, apparently, $\diam(T') = \max\{m_1, m_2, m_3, m_4\}$ and, thus,
  \[\diam(T') = \begin{cases}
    \ell_{2,3}, &\text{if } p = 0 \text{ and } s > 2 \text{ and } k_1-\odd(T_1) = k_3-\odd(T_3),\\
    \ell_{1,2}, &\text{otherwise},\\
  \end{cases}\]
  according to the considerations above.
  In the first case, where
  \[\diam(T') = \ell_{2,3} = k-\odd(k) + k_2 - \odd(T_2) + k_3 - \odd(T_3) - 2,\]
  the value $k - \odd(k) + k_2 + k_3 - 2$ is even because of $p=0$, and, therefore, $\odd(T') = \odd(T_2) + \odd(T_3) - 2\odd(T_2)\cdot\odd(T_3)$.
  Since $k_1 \geq k_2 \geq k_3$ are of the same parity and because
  \[k_1-\odd(T_1) \geq k_2-\odd(T_2) \geq k_3-\odd(T_3) = k_1-\odd(T_1),\]
  it must be the case that $k_1 = k_2 = k_3$ and $\odd(T_1) = \odd(T_2) = \odd(T_3)$.
  Moreover, Step~(\ref{const:case:connected}\ref{const:case:connected:k}) sets $k = k_1 + k_2 - 2\cdot\odd(T_1)\cdot\odd(T_2)$ in this case.
  Because $k_1 = k_3$ and $\odd(T_1) = \odd(T_3)$, we get
  \begin{align*}
    \ell_{2,3} &= k - 2 + (k_1 + k_2 - 2\cdot\odd(T_1)\cdot\odd(T_2)) - (\odd(T_2) + \odd(T_3) - 2\odd(T_2)\cdot\odd(T_3)) = 2k - 2 - \odd(T')
  \end{align*}
  and, hence, $\rad(T') = k-1$.
  Moreover, since $n \geq |V(G_1)| + |V(G_2)| + 1$ (because $u \in V(G) \setminus (V(G_1) \cup V(G_2))$), we have
  \begin{align*}
    (k-1) + (k-1) &\leq (|V(G_1)| + (k - k_1)) + (|V(G_2)| + (k - k_2))\\
    &= (|V(G_1)| + |V(G_2)|) + (k - k_1 - k_2) + k\\
    &\leq (n - 1) - 2\cdot\odd(T_1)\cdot\odd(T_2) + k\\
    &\leq n + (k-1),
  \end{align*}
  which means that $k-1 \leq n$.
  
  In the other case, where
  \[\diam(T') = \ell_{1,2} = k + k_1 - \odd(T_1) + k_2 - \odd(T_2) - 3,\]
  we distinguish between $p = 0$ and $p = 1$.
  If $p = 0$ then $k + k_1 + k_2 - 3$ is odd and $\odd(T') = 1 - \odd(T_1) - \odd(T_2) + 2\cdot\odd(T_1)\cdot\odd(T_2)$.
  Moreover, Step~(\ref{const:case:connected}\ref{const:case:connected:k}) sets $k = k_1 + k_2 - 2\cdot(\odd(T_1)+\odd(T_2)-\odd(T_1)\cdot\odd(T_2))$, and we get
  \begin{eqnarray*}
    \ell_{1,2} &=& k - 2 + (k_1 + k_2 - 2\cdot(\odd(T_1) + \odd(T_2) - \cdot\odd(T_1)\cdot\odd(T_2)) - (1 - \odd(T_1) - \odd(T_2) + 2\cdot\odd(T_1)\cdot\odd(T_2))\\
                &=& 2k - 2 - \odd(T'),
  \end{eqnarray*}
  which again means, $\rad(T') = k-1$.
  Otherwise, if $p = 1$ then $k + k_1 + k_2 - 3$ is even and $\odd(T') = \odd(T_1) + \odd(T_2) - 2\cdot\odd(T_1)\cdot\odd(T_2)$.
  In this final case, Step~(\ref{const:case:connected}\ref{const:case:connected:k}) defines $k = k_1 + k_2 - 1 - 2\cdot\odd(T_1)\cdot\odd(T_2))$ and, thus,
  \begin{eqnarray*}
    \ell_{1,2} &=& k - 2 + (k_1 + k_2 - 1 - 2\cdot\odd(T_1)\cdot\odd(T_2)) - (\odd(T_1) + \odd(T_2) - 2\cdot\odd(T_1)\cdot\odd(T_2))= 2k - 2 - \odd(T')
  \end{eqnarray*}
  and we can conclude that $\rad(T') = k-1$ holds in every case.
  Similarly to the previous case, we also get
  \begin{align*}
    (k-1) + (k-1) &\leq (|V(G_1)| + (k - k_1)) + (|V(G_2)| + (k - k_2))\\
    &= (|V(G_1)| + |V(G_2)|) + (k - k_1 - k_2) + k\\
    &\leq (n - 1) - (2\cdot\odd(T_1)\cdot\odd(T_2)+1) + k\\
    &\leq n + (k-1),
  \end{align*}
  and, again, $k-1 \leq n$.

  \medskip
  With $T'$, the algorithm has created a $k$-leaf root of $H$ with $\rad(T') = k-1 \leq n$ and $\odd(k) = p$.
  Step~(\ref{const:case:connected}\ref{const:case:connected:assembly}) extends $T'$ to a tree $T$ by attaching the universal vertex $u$ of $G$ as a leaf pendent to a center vertex of $T'$.  
  Obviously, $T$ is a $k$-leaf root of $G$ with $\odd(k) = p$ and property \textup{(T1)} $\rad(T) = k-1 \leq n$ because $\dist_T(x,y) = \dist_{T'}(x,y)$ for all $x,y \in V(H) = V(G) - u$ and
  \[\dist_T(u,x) \leq \rad(T) + 1 = \rad(T')+1 = (k-1)+1 = k.\]
  Moreover, $T$ has property \textup{(T2)} $\dmin_T = 1+\odd(T)$.
  If $\odd(T) = 0$ then the center of $T$ is a single vertex $z$, which is also the min-max center.
  Since $z$ is adjacent to the leaf $u$, we get $\dmin_T = 1$.
  Otherwise, if $\odd(T) = 1$ then $T$ has two center vertices $z_1$ and $z_2$ where, without loss of generality, $u$ is attached to $z_1$.
  That makes $z_2$ the min-max center with distance $\dmin_T = 2$ to leaf $u$.
  This concludes the proof of Proposition~\ref{lem:correct}.

  \medskip  
  With the following proposition, we summarize what we have so far and cover the disconnected case.
  \begin{proposition}\label{thm:correct}
    Let $G$ be a {\TPG} on $n$ vertices and without true twins and let $p \in \{0,1\}$ be a given parity.
    Then $(T,k) = \rho(G,p)$ provides a $k$-leaf root $T$ of $G$ with $\odd(k) = p$ and $k \leq n+1$.
    If $G$ is connected then $(T,k)$ satisfies \textup{(T1)} and \textup{(T2)}.
  \end{proposition}
  After the Propositions~\ref{obs:star} and~\ref{lem:correct}, it remains to handle disconnected input.
  In this case, $G$ has $s \ge 0$ non-trivial connected components $G_1, \ldots, G_s$ and $t \ge 0$ isolated vertices $v_1,\ldots,v_t$ such that $s+t \ge 2$.

  If $s = 0$ then $G$ consists of $n \geq 2$ isolated vertices and, thus, Step~(\ref{const:case:disconnected}\ref{const:case:disconnected:widen}) sets $k = 2$ for even parity and $k = 3$, otherwise.
  Then, by Lemma~\ref{lem:disconnected}, Step~(\ref{const:case:disconnected}\ref{const:case:disconnected:assembly}) returns a $k$-leaf root $T$ of $G$.

  For $s > 0$, every $G_i$ is a connected {\TPG} without true twins.
  By Proposition~\ref{lem:correct}, Step~(\ref{const:case:disconnected}\ref{const:case:disconnected:recursion}) provides a $k_i$-leaf root $T_i$ of parity $p$ for every $G_i$, $1\le i\le s$.
  The value assigned to $k$ in Step~(\ref{const:case:disconnected}\ref{const:case:disconnected:widen}) is the maximum $k_m = \max\{k_1, \dots, k_s\}$.
  Since $k_m \leq |V(G_m)| + 1$, we have $k \leq m + 1$.
  According to Lemma~\ref{lem:subdividing}, the same step produces a $k$-leaf root $T'_i = \eta\left(T_i, \tfrac{k-k_i}{2}\right)$ of $G_i$ for all $i \in \{1, \dots, s\}$.
  Finally, by Lemma~\ref{lem:disconnected}, Step~(\ref{const:case:disconnected}\ref{const:case:disconnected:assembly}) returns a $k$-leaf root $T$ of $G$.
  This proves Proposition~\ref{thm:correct}.
  
  \medskip
  The following, last proposition finalized the proof of Theorem~\ref{cor:correct_optimal} because it shows that the $k$-leaf root $T$ returned by $\rho(G, p)$ is optimal with respect to both, the value of $k$ and the diameter of $T$.
  \begin{proposition}\label{thm:optimal}
    Let $G$ be a {\TPG} without true twins and let $p \in \{0,1\}$ be a given parity.
    Then $(T,k) = \rho(G,p)$ provides a $k$-leaf root $T$ of $G$ that is \emph{optimal} for parity $p$ (hence, $\odd(k) = p$).
    If $G$ is connected then $(T,k)$ satisfies \textup{(T3)}.
  \end{proposition}

  That $(T,k) = \rho(G,p)$ provides a $k$-leaf root $T$ of $G$ with $\odd(k) = p$ has just been established.
  To prove Proposition~\ref{thm:optimal}, we only need to show that $T$ fulfills \textup{(T3)} for connected input graphs $G$ and that $k$ is always parity-optimal.
  
  We begin with the proof of \textup{(T3)} and, like for Proposition~\ref{thm:correct}, this works by induction on the vertex number of $G$. 
  Again, for the base case, Proposition~\ref{obs:star} handles the star with two leaves, the smallest connected {\TPG} without true twins and more than one vertex.
  Proposition~\ref{obs:star} settles the theorem for all stars.
    
  \medskip
  Now, let $G$ have more than three vertices.
  That $G$ is connected and not a star leads us to Case~\ref{const:case:connected} of $\rho$, again.
  Borrowing the argumentation from the proof of Proposition~\ref{thm:correct}, we immediately get that
  \begin{itemize}
    \item
    $G$ has a unique universal cutvertex $u$ and $H = G-u$ has $s \ge 1$ non-trivial connected components $G_1, \ldots, G_s$, each a {\TPG} without true twins, and $t \ge 0$ isolated vertices $v_1,\ldots,v_t$ such that $s+t \ge 2$,
    \item
    by induction hypothesis, Step~(\ref{const:case:connected}\ref{const:case:connected:recursion}) provides $k_i$-leaf roots $T_i$ of parity $p$ for all $G_i, 1 \leq i \leq s$, which are optimal and satisfy \textup{(T3)}, here,
    \item
    $k_1 \geq k_2 \geq \dots \geq k_s$, if, without loss of generality, we assume
    \[\diam(T_1) \geq \diam(T_2) \geq \dots \geq \diam(T_s),\]
    \item
    Step~(\ref{const:case:connected}\ref{const:case:connected:widen}) produces $k$-leaf roots $T'_i = \eta\left(T_i, \tfrac{k-k_i}{2}\right)$ of all $G_i, 1 \leq i \leq s$ with
    \begin{align*}
      k &= \begin{cases}
      k_1 + 2(1 - \odd(T_1)),                                         &\text{if } s = 1,\\
      k_1 + k_2 - 1 - 2\cdot\odd(T_1)\cdot\odd(T_2)                   &\text{if } s \geq 2 \text{ and } p = 1,\\
      k_1 + k_2 - 2\cdot\odd(T_1)\cdot\odd(T_2),                      &\text{if } s > 2, p = 0 \text{ and } k_1-\odd(T_1) = k_3-\odd(T_3),\\
      k_1 + k_2 - 2 \cdot (\odd(T_1)+\odd(T_2) - \odd(T_1)\cdot\odd(T_2)),&\text{otherwise,}
    \end{cases}\\
      \diam(T'_i) &= \diam(T_i)+k-k_i, \text{ and}\\
      \rad(T'_i) &= \rad(T_i) + \tfrac{k-k_i}{2} = \tfrac{k+k_i}{2} - 1,\\
  \end{align*}
    \item
    Step~(\ref{const:case:connected}\ref{const:case:connected:widen}) produces $T' = \mu(k, T'_1, \dots, T'_s, v_1, \dots, v_t)$ with
    \[\diam(T') = \begin{cases}
      k + k_1 - \odd(T_1) - 1, &\text{if } s = 1,\\
      k + k_2 - \odd(T_2) + k_3 - \odd(T_3) - 2, &\text{if } s > 2 \text{ and } p = 0 \text{ and } k_1-\odd(T_1) = k_3-\odd(T_3),\\
      k + k_1 - \odd(T_1) + k_2 - \odd(T_2) - 3, &\text{otherwise, and}
    \end{cases}\]
    \item
    $(T,k) = \rho(G, p)$, as assembled in Step~(\ref{const:case:connected}\ref{const:case:connected:assembly}), represents a $k$-leaf root of $G$ with $\diam(T) = \diam(T')$.
  \end{itemize}
  
  Next, consider an arbitrary $k'$-leaf root $R$ of $G$ and, for all $i \in \{1, \dots, s\}$, let $R_i$ be the smallest subtree of $R$ with leaf set $V(G_i)$.
  Observe that the $R_i$-trees are pairwise vertex-disjoint.
  Then $R_i$ is a $k'$-leaf root of $G_i$ and, since $T_i$ has \textup{(T3)},
  \[\diam(R_i) \ge \diam(T_i) + k'-k_i = (\diam(T_i)-k+k_i) + k'-k_i = \diam(T'_i) + k'-k.\]
    
  Let $u_i$ be the unique universal vertex of $G_i$ by Proposition~\ref{lem:wolk_plus}.
  By Lemma~\ref{lem:A}, there is a center vertex $z_i$ of $R_i$ with $\dist_{R_i}(u_i,z_i) \le k' - \rad(R_i)$.
  Let $x_i$ and $y_i$ be the two leaves of a diametral path in $R_i$.
  Then, the center vertex $z_i$ of $T_i$ is on the $x_i,y_i$-path and $\dist_{R_i}(u_i,z_i) = \dist_R(u_i,z_i)$.

  \medskip
  Suppose $s=1$.
  Then $\dist_R(z_1,v_1) > \rad(R_1)$ for, otherwise,
  \[\dist_R(u_1,v_1) \le \dist_R(u_1,z_1) + \dist_R(z_1,v_1) \le (k'-\rad(R_1)) + \rad(R_1) = k',\]
  and $u_1$ and $v_1$ would be adjacent in $G$.
  In~$R$, the longest path between $v_1$ and a vertex in $\{x_1,y_1\}$ must contain~$z_1$, and it has the length $\ell$ that fulfills
  \[\ell \ge \dist_R(v_1,z_1) + \rad(R_1) > 2\rad(R_1) \ge \diam(R_1) \geq \diam(T'_1) + k'-k.\]
  Hence, $\diam(R) > \diam(T'_1) + k'-k$ and we conclude
  \begin{align*}
    \diam(R) &\geq \diam(T'_1) + k'-k + 1 = (2\rad(T'_1)-\odd(T'_1)) + k'-k + 1\\
              &= (k+k_1-2) - \odd(T'_1) + k'-k + 1\\
              &= (k + k_1 - \odd(T_1) - 1)  + k'-k\\
              &= \diam(T) + k'- k,
  \end{align*}
  as claimed.
  
  \medskip
  If $s \ge 2$ then $\dist_R(z_i,z_j) \ge \rad(R_i) + \rad(R_j) - k'+1$ for all $1 \le i < j \le s$ since, otherwise,
  \begin{align*}
    \dist_R(u_i,u_j) &\le \dist_R(u_i,z_i) + \dist_R(z_i,z_j) + \dist_R(z_j,u_j)\\
                      &\le (k'-\rad(R_i)) + (\rad(R_i)+\rad(R_j)-k') + (k'-\rad(R_j))= k',
  \end{align*}
  and~$u_i$ and~$u_j$ would be adjacent in~$G$.   
  In~$R$, the longest path going between a vertex in $\{x_i,y_i\}$ and a vertex in $\{x_j,y_j\}$ contains~$z_i$ and~$z_j$. 
  The length $\ell_{i,j}$ of this path fulfills 
  \begin{align*}
  \ell_{i,j} &\ge \rad(R_i) + \dist_R(z_i,z_j) + \rad(R_j)\\ 
              &\ge \rad(R_i) + \big(\rad(R_i)+\rad(R_j)-k+1\big) + \rad(R_j)\\ 
              &= \diam(R_i) + \odd(R_i) + \diam(R_j) + \odd(R_j) - k'+1.
  \end{align*}
  Because $\diam(R) \geq \ell_{1,2}$ and $\odd(R_1) \geq 0$ and $\odd(R_2) \geq 0$, we find that
  \begin{align*}
    \diam(R) &\ge \diam(R_1) + \diam(R_2) - k'+1\\
              &\ge (\diam(T_1) + k'-k_1) + (\diam(T_2) + k'-k_2) - k'+1\\
              &= (2\rad(T_1) - \odd(T_1) - k_1) + (2\rad(T_2) - \odd(T_2) - k_2) + k' + 1\\
              &= (2k_1-2 - \odd(T_1) - k_1) + (2k_2-2 - \odd(T_2) - k_2) + k' + (k-k) + 1\\
              &= (k + k_1 - \odd(T_1) + k_2 - \odd(T_2) - 3) + k' - k.
  \end{align*}
  Hence, if $p = 1$ or $s = 2$ or ($s > 2$ and) $k_1-\odd(T_1) > k_3-\odd(T_3)$ then (as $\diam(T) = k + k_1 - \odd(T_1) + k_2 - \odd(T_2) - 3$, here) we already get
  \[\diam(R) \geq \diam(T) + k' - k.\]
  
  In the remaining case, we have $p = 0$ and $s > 2$ and $k_1-\odd(T_1) = k_3-\odd(T_3)$.
  Like in the proof of Theorem~\ref{thm:correct}, we can conclude that $k_1 = k_2 = k_3$ and $\odd(T_1) = \odd(T_2) = \odd(T_3)$, here.
  Thus, $\diam(T_1) = \diam(T_2) = \diam(T_3)$ and $\diam(T'_1) = \diam(T'_2) = \diam(T'_3)$, too.
  If we write $\delta = \diam(T'_1) + k'-k$ then $\diam(R_i) \geq \delta$ for all $i \in \{1,2,3\}$.

  Firstly, we show that at least one of the three length $\ell_{1,2}, \ell_{1,3},\ell_{2,3}$ is at least $2\delta - k' + 2$.
  If $\diam(R_i)+\odd(R_i) > \delta$ for some $i\in\{1,2,3\}$, then 
  \begin{align*}
    \ell_{i,j} &\ge \diam(R_i) + \odd(R_i) + \diam(R_j) + \odd(R_j) - k'+1\\
                &\ge (\delta + 1) + \delta - k'+1\\
                &\ge 2\delta - k' + 2
  \end{align*}
  for all $j \in \{1,2,3\} \setminus \{i\}$ and we are done, already.
  Therefore, consider that $\diam(R_i)+\odd(R_i) \le \delta$ for all $i \in \{1,2,3\}$.
  Then, since and $\delta \leq \diam(R_i) \leq \delta-\odd(R_i)$ for all $i \in \{1,2,3\}$, we have
  \[\diam(R_1)+\odd(R_1) = \diam(R_2)+\odd(R_2) = \diam(R_3)+\odd(R_3) = \delta.\]
  Hence, $\rad(R_1)=\rad(R_2)=\rad(R_3) = \tfrac{\delta}{2}$, which implies that $\delta$ is even.
  If $\dist_R(z_i,z_j) \ge \delta - k'+2$ for any $i,j \in \{1,2,3\}$, then
  \begin{align*}
    \ell_{i,j} &\ge \rad(R_i) + \dist_R(z_i,z_j) + \rad(R_j)\\
              &\ge \tfrac{\delta}{2} + (\delta - k'+2) + \tfrac{\delta}{2}\\
              &= 2\delta - k'+2
  \end{align*}
  and done, again.
  So, finally, assume that $\dist_R(z_i,z_j) \le \delta-k'+1$ for all $i \in \{1,2,3\}$.
  Then, since we have
  \[\dist_R(z_i,z_j) \ge \rad(R_i) + \rad(R_j) - k'+1 = \tfrac{\delta}{2} + \tfrac{\delta}{2} - k'+1\]
  from before, $\dist_R(z_1,z_2) = \dist_R(z_1,z_3)=\dist_R(z_2,z_3) = \delta - k'+1$.
  Let $w$ be the last vertex $w$ on the $z_1,z_2$-path in $R$ that is closest to $z_3$ and let $\delta_i$ be the distances between $z_i$ and $w$ for all $i \in \{1,2,3\}$.
  Then, for all $1\le i<j\le 3$, it must be
  \[\delta_i + \delta_j = \dist_R(z_i,z_j) = \delta-k'+1.\]
  Solving the system of the three linear equations above yields
  \[\delta_1 = \delta_2 = \delta_3 = \tfrac{1}{2}(\delta-k'+1),\]
  which implies that $k'$ is odd, a contradiction to $\odd(k') = \odd(k) = p = 0$.

  So at least one length, say $\ell_{1,2}$, is at least $2\delta - k'+2$.
  This means that
  \begin{align*}
    \diam(R) &\geq \ell_{1,2} \geq 2\delta - k'+2\\
              &= (\diam(T'_1) + k'-k) + (\diam(T'_3) + k'-k) - k'+2\\
              &= (\diam(T_1) + k-k_1 + k'-k) + (\diam(T_3) + k - k_3 + k'-k) - k'+2\\
              &= (2\rad(T_1) - \odd(T_1) + k'-k_1) + (2\rad(T_3) - \odd(T_3) + k'-k_3) - k'+2\\
              &= (2\rad(T_1) - \odd(T_1) -k_1) + (2\rad(T_3) - \odd(T_3) -k_3) + k' + (k-k) + 2\\
              &= k + (k_1 - 2 - \odd(T_1)) + (k_3 - 2 - \odd(T_3)) + k'-k + 2\\
              &= (k + k_1 - \odd(T_1) + k_3 - \odd(T_3) - 2) + k' - k.
  \end{align*}
  Because here, $\diam(T) = k + k_1 - \odd(T_1) + k_3 - \odd(T_3) - 2$, we have now shown that $\diam(R) \geq \diam(T) + k' - k$.

  \medskip
  It remains to prove that the value of $k$ in $(T,k) = \rho(G,p)$ is optimal with respect to leaf roots of $G$ with parity $p$.
  For the connected case, this can be seen as follows.
  By Corollary~\ref{cor:A}, $2k'-2 \ge \diam(R)$ holds for every $k'$-leaf roots $R$ of~$G$. 
  Hence,
  \[2k'-2 \ge \diam(R) \ge \diam(T)+k'-k \ge 2k-3+k'-k = k+k'-3,\]
  since $\diam(T) = 2\rad(T)-\odd(T) \geq 2(k-1) - 1 = 2k-3$.
  Now, $2k'-2 \geq k+k'-3$ implies that $k'\ge k-1$ and, because both, $k$ and $k'$, have the same parity, $k'\ge k$. 
  That is, $T$ is an optimal $k$-leaf root of~$G$ with parity $p$.

  \medskip
  Finally, let $G$ be a disconnected graph with $s \ge 0$ non-trivial connected components $G_1, \dots, G_s$ and $t \ge 0$ isolated vertices $v_1, \dots, v_t$ such that $s+t \geq 2$.
  By Theorem~\ref{thm:correct}, $(T,k) = \rho(G, k)$ results from finding $(T_1, k_1) = \rho(G_1, p), \dots, (T_s, k_s) = \rho(G_s, p)$ in Step~(\ref{const:case:disconnected}\ref{const:case:disconnected:recursion}) and assembling them to a $k$-leaf root $T$ for $G$ with $k = \max\{k_1, \dots, k_s, p+2\}$ in the Steps~(\ref{const:case:disconnected}\ref{const:case:disconnected:widen}) and~(\ref{const:case:disconnected}\ref{const:case:disconnected:assembly}).
  After handling the connected case above, we know that $T_i$ is an optimal $k_i$-leaf root of $G_i$ with parity $p$ for all $i \in \{1, \dots, s\}$.
  
  Clearly, if $s = 0$ then a $(p+2)$-leaf root is the best possible for parity $p$, and we are done.
  Otherwise, assume, without loss of generality, that $k = k_1$.
  Moreover, let $R$ be a $p$-parity $k'$-leaf root of $G$ and, for all $i \in \{1, \dots, s\}$, let $R_i$ be the smallest subtree of $R$ with leaf set $V(G_i)$.
  Then $R_1$ is a $k'$-leaf root of $G_1$.
  Since $k_1$ is optimal for $G_1$, we have that $k' \geq k_1 = k$, thus, $T$ is an optimal $k$-leaf root of~$G$ with parity $p$.
  This proves Proposition~\ref{thm:optimal} and, with it, also Theorem~\ref{cor:correct_optimal}.
\qed\end{proof}
  

\repetition{Lemma}{lem:subdividing:complexity}{
  Let $T$ be a compressed tree with $n$ leaves and with explicitly given min-max center $z$, center $Z$, diameter $\diam(T)$, and leaf-distance $\dmin_T$.
  For all integers $\delta \geq 0$, the compressed tree $T' = \eta(T,\delta)$ can be computed in $\cO(n)$ time together with min-max center $z'$, center $Z'$, diameter $\diam(T')$, and leaf distance $\dmin_{T'}$ of $T'$.
}
\begin{proof}
  By definition, $T'$ is obtained from $T$ by replacing all $n$ pendant edges with new paths of length $\delta+1$, each.
  In the compressed encoding, this takes just $n$ modifications of the weight of edges.
  More precisely, if $v(\ell)x$ is a pendent edge of $T$ representing a path of length $\ell$ that ends at leaf $x$ (where we simply take $\ell=1$ for an unweighted edge) then $T'$ has the same pendent edge $v(\ell')x$ with the new weight $\ell' = \ell + \delta$.
  The $n$ constant-time modifications take $\cO(n)$ time, altogether.
  
  According to Lemma~\ref{lem:subdividing}, $z' = z$, $Z' = Z$, $\diam(T') = \diam(T)+2\delta$ and $\dmin_{T'}=\dmin_{T}+\delta$, which takes just constant time to compute.
\qed\end{proof}
  
\repetition{Lemma}{lem:disconnected:complexity}{
  Let $s \geq 0$ be an integer and, for all $i \in \{1, \dots, s\}$, let $T_i$ be a given, compressed tree with explicitly given min-max center $z_i$, center $Z_i$, diameter $\diam(T_i)$, and leaf-distance $\dmin_{T_i}$.
  For all integers $k \geq 2$ and vertices $v_1, \dots, v_t$, the merged compressed tree $T' = \mu(k, T_1, \dots, T_s, v_1, \dots, v_t)$ with center $Z'$ and diameter $\diam(T')$ can be computed in $\cO(s + t)$ time.
}
\begin{proof}
  According to definition, the computation of $T'$ by $\mu$ starts with finding the critical index $m$, which is the smallest number in $\{1, \dots, s\}$ with $\dmin_{T_m} = \min\{\dmin_{T_i} \mid 1 \leq i \leq s\}$.
  Since $\dmin_{T_i}$ is explicitly known for all given trees, $m$ is found after iterating the input once in $\cO(s)$ time.
  
  Then, every tree $T_i, 1 \leq i \leq s$ is attached at $z_i$ to a new vertex $c$ by a path of at most $\tfrac{k-\odd(k)}{2}+1-\dmin_{T_i}$ edges and every vertex $v_j, 1 \leq j \leq t$ is attached to $c$ with a path of length $\tfrac{k-\odd(k)}{2}+1$.
  Each of these paths is represented by one edge of the respective weight and, hence, this takes $\cO(s + t)$ time.
  
  \medskip
  It remains to show that $Z'$ and $\diam(T')$ can be computed alongside $T'$ without consuming essentially more computing time.
  Both parameters require the identification of a diametral path in $T'$.
  We recall that every diametral path $P$ connects two leaves of the tree $T'$ and, therefore, is either entirely included in one of the trees $T_1, \dots, T_s$ or contains $c$ and connects leaves that stem from different input trees or vertices $v_1, \dots, v_t$.
  To be able to cover the first case, we search a tree $T_a$ with $\diam(T_a) = \min\{\diam(T_i) \mid 1 \leq i \leq s\}$.
  This works in $\cO(s)$ time since the diameters are given.
  If, at the end, $P$ turns out to be part of $T_a$ then $Z' = Z_a$ and $\diam(T') = \diam(T_a)$, which is computed in $\cO(1)$ time.
  
  For the other possibility, we firstly think of the case where $P$ connects two leaves $v_i$ and $v_j$, $1 \leq i < j \leq t$.
  Then, obviously, $\diam(T') = k-\odd(k)+2 \geq \diam(T_a)$ and $Z' = \{c\}$, which can be decided and computed in $\cO(1)$ time.
  
  Secondly, one or both end vertices of $P$ may be leaves of the trees $T_1, \dots, T_s$.
  To efficiently find $P$ under this condition, we hook additional computations into the iteration of these trees during the evaluation of the $\mu$-operation.
  More precisely, while a tree $T_i$ is attached to $c$ by a path of length $p_i$ (that is, $p_i = \tfrac{k+\odd(k)}{2} - \dmin_{T_i}$ for $i = m$ and, otherwise, $p_i = \tfrac{k-\odd(k)}{2}+1-\dmin_{T_i}$), we memorize the $T'$-distance
  \[d_i = \max\{\dist_{T'}(c,v) \mid v \text{ is leaf in } T_i\} = p_i + \rad(T_i) = p_i + \left\lceil\tfrac{\diam(T_i)}{2}\right\rceil\]
  between $c$ and a farthest leaf of $T_i$.
  As $\diam(T_i)$ is known in advance for all $i \in \{1, \dots, s\}$, this takes just $\cO(1)$ additional time per iteration and, thus, $\cO(s)$ time in total.
  We also define
  \[d_{0} = \begin{cases}
    \dist_{T'}(c,v_1) = \tfrac{k-\odd(k)}{2}+1, &\text{if } t > 0,\\
    0, &\text{otherwise},
  \end{cases}\]
  the $T'$-distance between $c$ and any present leaf $v_1, \dots, v_t$.
  Observe for later that, if $t > 0$ then
  \begin{align*}
    d_i - d_0 &= p_i + \rad(T_i) - \tfrac{k-\odd(k)}{2}+1\\
      &\geq \tfrac{k+\odd(k)}{2} - \dmin_{T_i} + \rad(T_i) - \tfrac{k-\odd(k)}{2}+1\\
      &\geq \tfrac{(k+\odd(k))-(k-\odd(k))}{2} + 1 \text{ since } \dmin_{T_i} \leq \rad(T_i)\\
      &= \odd(k) + 1 \geq 1,
  \end{align*}
  hence, $d_i > d_0$ for all $i \in \{1, \dots, s\}$.

  After the evaluation of $\mu$ and the associated completion of $T'$, we determine values $i_1$ and $i_2$ such that
  \[d_{i_1} = \max\{d_i \mid 0 \leq i \leq s\} \text{ and } d_{i_2} = \max\{d_i \mid 0 \leq i \leq s, i \not = i_1\},\]
  which takes $\cO(s)$ time.
  If a diametral path $P$ of $T'$ runs through $c$ and has at least one end vertex in one of the trees $T_1, \dots, T_s$ then, obviously, the length of this path is $d = d_{i_1} + d_{i_2}$.
  We can detect this situation in $\cO(1)$ time by checking $d > \diam(T_a)$ and $t>1 \Rightarrow d > k-\odd(k)+2$.
  In this case, we infer that $\diam(T') = d$ and $Z'$ is on $P$.
  More precisely, if $Z' = \{z'_1, z'_2\}$ (with $z'_1 = z'_2$ if $\odd(T') = 1$) then $Z'$ is situated on the longer subpath of $P$ with $z'_1$ at distance $\delta_1 = \left\lceil\tfrac{d_{i_1} - d_{i_2}}{2}\right\rceil$ from $c$ and $z'_2$ at distance $\delta_2 = \left\lfloor\tfrac{d_{i_1} - d_{i_2}}{2}\right\rfloor$.
  
  To conclude the proof, we show that $Z'$ is situated on the recently inserted weighted edge $c(p_{i_1})z_{i_1}$, or more precisely, $p_{i_1} \geq \delta_1 \geq \delta_2 \geq 0$.
  That $\delta_2 \geq 0$ holds since, otherwise, $d_{i_1} < d_{i_2}$. 
  Similarly, we note that, if one of the end vertices of $P$ is a leaf of $v_1, \dots, v_t$ then $i_2 = 0$ since, otherwise, $d_0 = d_{i_1} < d_{i_2}$.
  So, $i_1$ is in $\{1, \dots, s\}$, and we further observe that
  \[2\rad(T_{i_1}) - \odd(T_{i_1}) = \diam(T_{i_1}) \leq \diam(T_a) < d_{i_1} + d_{i_2} = p_{i_1} + \rad(T_{i_1}) + d_{i_2},\]
  because $T_a$ has the largest diameter among the given $s$ trees.
  This means that $\rad(T_{i_1}) \leq p_{i_1} + d_{i_2}$.
  Using this for the estimation of $\delta_1$, we get
  \begin{align*}
    \delta_1 &= \left\lceil\tfrac{d_{i_1} - d_{i_2}}{2}\right\rceil
    \leq \tfrac{p_{i_1} + \rad(T_{i_1}) - d_{i_2}}{2}
    \leq \tfrac{p_{i_1} + (p_{i_1} + d_{i_2}) - d_{i_2}}{2}\\
    &= p_{i_1}.
  \end{align*}
  Hence, if $\delta_1 = \delta_2 = 0$ then $z'_1 = z'_2 = c$.
  If $\delta_1 = p_{i_1}$ then $z'_1 = z_{i_1}$ and, if $\delta_1 > \delta_2$ in this case, then the edge $c(p_{i_1})z_{i_1}$ is split into the weighted edge $c(p_{i_1}-1)z'_2$ and the edge $z'_2z'_1$.
  Otherwise, if also $\odd(T') = 0$, then the center $z'_1 = z'_2$ is inserted by splitting the edge $c(p_{i_1})z_{i_1}$ into the weighted edges $c(\delta_1)z'_1$ and $z'_1(p_{i_1}-\delta_1)z_{i_1}$.
  In the final case, $c(p_{i_1})z_{i_1}$ is replaced by inserting the center edge $z'_1z'_2$ between $c(\delta_2)z'_2$ and $z'_1(p_{i_1}-\delta_1)z_{i_1}$.
  Altogether, this implies that $Z'$ and $\diam(T')$ are always found in $\cO(s)$ time.
\qed\end{proof}
  
\repetition{Theorem}{thm:main:odd_even}{
  Given a chordal cograph $G$ on $n$ vertices and $m$ edges and a parity $p$ (either odd or even), a (compressed) $\kappa$-leaf root of $G$ with minimum integer $\kappa$ of parity $p$ can be computed in $\cO(n + m)$ time.
}
\begin{proof}
  Let $G$ be a {\TPG} and $p \in \{0,1\}$ a parity.
  We can assume that $G$ is free of true twins as, otherwise, we can remove $x$ for every pair $x,y$ of true twins in $G$ in linear time and, when $(T, k)$ has been computed, insert every $x$ pendant to the parent of the corresponding twin leaf $y$ into $T$.
  The cotree $\cT$ of $G$ is computed in Line~\ref{algo:init_ready}, which works on linear time $\cO(n + m)$ as discussed in Section~\ref{sec:basic}.

  We begin with connected input graphs and use complete induction on $n$ to show that, then, Algorithm~\ref{alg:opt} puts $(T,k) = \rho(G,p)$ on top of the stack $\cS$ before finishing the last visited node of $\cT$ by the post-order-loop in Lines~\ref{algo:loop_start} to~\ref{algo:loop_end}, which requires $\cO(n)$ steps.
  Since $G$ is connected and free of true twins, the induction start is at $n = 3$ with $G$ the star with center $u$ and two leaves $v_1, v_2$.
  For any star with leaves $v_1, \dots, v_t$, the cotree has a $\otimes$-root $X$ with leaf child $u$ and a $\oplus$-child that has the leaf-children $v_1, \dots, v_t$, only.
  This means, Algorithm~\ref{alg:opt} iterates just one $\otimes$-node (the root $X$), heads directly into Line~\ref{algo:start_base_case}, and, accordingly, handles this case like Case~\ref{const:base} of the $\rho$-operation.
  Afterwards, it pushes $(T,k) = \rho(G,p)$ on top of the stack $\cS$ before the for-loop finishes $X$, the last visited node of $\cT$.
  It is easy to see that, the $(4-p)$-leaf root $T$ on $t+1$ vertices is constructed in $\cO(t) = \cO(n)$ time.

  \medskip
  Next for the induction step, where $n > 3$ and $G$ is not a star.
  The algorithm enters Line~\ref{algo:start_connected_case} because $G$ is connected and, thus, the cotree $\cT$ has a $\otimes$-root $X$ with leaf-child $u$ and a $\oplus$-child $Y$ with, in turn, $\otimes$-children $Z_1, \dots, Z_s$, $s \geq 1$.
  Due to post-order traversal, the algorithm consecutively iterates through the subtrees $\cT_{Z_1}, \dots, \cT_{Z_s}$ before entering $Y$ and then $X$.
  Since the graphs $G_1, \dots, G_s$ that correspond to $\cT_{Z_1}, \dots, \cT_{Z_s}$ are connected and have strictly fewer vertices than $G$ (they all miss $u$), the induction assumption states that $\cS$ carries the sequence $(T_1, k_1) = \rho(G_1,p), \dots, (T_s, k_s) = \rho(G_s,p)$ on its top.
  Hence, Algorithm~\ref{alg:opt} retrieves this list of interim results in Line~\ref{algo:start_connected_case}.
  Although unmentioned for clarity in the listing of Algorithm~\ref{alg:opt}, every tree $T_i$, $1 \leq i \leq s$ is retrieved together with the explicit information on center, min-max center, diameter, and leaf distance.
  Altogether, this takes $\cO(|V(G_1)|) + \dots + \cO(|V(G_s)|) = \cO(n)$ time.

  After that, it is easy to check that $k$ is computed exactly like in Case~\ref{const:case:connected} in the listing of the $\rho$-operation in Section~\ref{sec:tplr}.
  Hence, Line~\ref{algo:connected_widening} computes the same extended $k$-leaf roots $T'_1, \dots, T'_s$ as $\rho$ and, by Lemma~\ref{lem:subdividing:complexity}, this takes $\cO(|V(G_1)|) + \dots + \cO(|V(G_s)|) = \cO(n)$ time.
  Notice that the subroutine $\eta$ includes min-max center $z'_i$, center $Z'_i$, diameter $\diam(T'_i)$, and leaf distance $\dmin_{T'_i}$ for every tree $T'_i$, $1 \leq i \leq s$.
  Based on that, the Lines~\ref{algo:connected_merging} and~\ref{algo:connected_attach} produce the same merged $k$-leaf root $T$ as $\rho$ does.
  According to Lemma~\ref{lem:disconnected:complexity} and because $s+t \leq n$, this needs just $\cO(s+t) = \cO(n)$ time.
  The subroutine $\mu$ explicitly returns only the center $Z = \{z_1, z_2\}$ and the diameter $\diam(T)$ with $T$.
  However, the min-max center $z$ and the leaf distance $\dmin_T$ are immediately determined in the next Line~\ref{algo:connected_attach} due to attaching $u$.
  If $\odd(T) = 0$ then $z = z_1 = z_2$ and $\dmin_T = 1$ and, otherwise, $z$ the center that is not chooses as the parent of $u$ and $\dmin_T = 2$.
  This means that $(T,k) = \rho(G,p)$, including $z$, $Z$, $\diam(T)$, and $\dmin_T$, is subsequently pushed onto $\cS$ in Line~\ref{algo:connected_push}.
  This finishes the processing of $X$, the last visited node of $\cT$, after $\cO(n)$ time in total.
  
  \medskip
  For a disconnected graph $G$, we also show that the loop in Lines~\ref{algo:loop_start} to~\ref{algo:loop_end} puts $(T,k) = \rho(G,p)$ on top of the stack $\cS$.
  The cotree $\cT$ of $G$ has a $\oplus$-root $X$, which has $\otimes$-children $Z_1, \dots, Z_s$, $s \geq 0$ and leaf-children $V_1, \dots, V_t$ such that $s+t \geq 2$.
  Notice that $X$, as the root of $\cT$, has no parent.
  Then, as argued before, the post-order traversal consecutively iterates through $\cT_{Z_1}, \dots, \cT_{Z_s}$ before entering the node $X$ in Line~\ref{algo:disconnected_start}.
  Since the graphs $G_1, \dots, G_s$ are connected, we have seen above that iterating the children of $X$ takes $\cO(|V(G_1)|) + \dots + \cO(|V(G_s)|) = \cO(n)$ time and puts the sequence $(T_1, k_1) = \rho(G_1,p), \dots, (T_s, k_s) = \rho(G_s,p)$ on top of the stack $\cS$.
  Again unmentioned, every tree $T_i$, $1 \leq i \leq s$ comes with explicitly designated min-max center, center, diameter and leaf distance.
  Altogether, this takes $\cO(|V(G_1)|) + \dots + \cO(|V(G_s)|) = \cO(n)$ time.
  Like in Case~\ref{const:case:disconnected} of the $\rho$-operation, the value of $k$ is the maximum among $\{k_1, \dots, k_s, p+2\}$.
  So, by Lemmas~\ref{lem:subdividing:complexity} and~\ref{lem:disconnected:complexity}, $(T,k) = \rho(G,p)$ is correctly computed in Lines~\ref{algo:disconnected_widening} and~\ref{algo:disconnected_merging} in $\cO(n)$ time.
  Then Line~\ref{algo:disconnected_push} puts $(T,k)$ on top of $\cS$, which finishes $X$, the last visited node of $\cT$.
  
  Since $\cO(n) \leq \cO(m)$, the Algorithm takes $\cO(n + m)$ time in each of the above cases to finally reach Line~\ref{algo:interim_result}, where $(T,k) = \rho(G, p)$ waits accessible at the top of $\cS$.
  Hence, because $G$ is a {\TPG} without true twins, Theorem~\ref{cor:correct_optimal} tells us that Algorithm~\ref{alg:opt} provides an optimal $k$-leaf root $T$ of $G$ with $\odd(k) = p$ at this point.
  This proves Theorem~\ref{thm:main:odd_even}.
\qed\end{proof}


\repetition{Lemma}{lem:family_roots}{
  For all integers $i \geq 1$, the graph $F_i$ is a (odd) $k_i$-leaf power for $k_i = 2^{i+2}-1$ but not a $(k_i-2)$-leaf power and a (even) $k'_i$-leaf power for $k'_i = k_i + 2^{i} - 1$ but not a $(k'_i-2)$-leaf power.
}
\begin{proof}
  We show by induction on $i$ that $(T_i, k_i) = \rho(F_i, 1)$ and $(T'_i, k'_i) = \rho(F_i, 0)$ with $\odd(T_i) = \odd(T'_i) = 0$.
  This proves the lemma because, by Theorem~\ref{cor:correct_optimal}, $\rho$ provides an optimal leaf root of the given parity.

  For the induction start at $i=1$, we leave it to the reader to check that $(T_1, k_1) = \rho(F_1, 1)$ and $(T'_1, k'_1) = \rho(F_1, 1)$ provide a $k_1$-leaf root $T_1$ of $F_1$ with $k_1 = 7$ and $\diam(T_1) = 8$ and a $k'_1$-leaf root $T'_1$ with $k'_1 = 8$ and $\diam(T'_1) = 10$.
  
  \medskip
  For the induction step let $i > 0$ and assume that $(T_{i-1}, k_{i-1}) = \rho(F_{i-1}, 1)$ with $\odd(T_{i-1}) = 0$ and $(T'_{i-1},k'_{i-1}) = \rho(F_{i-1}, 0)$ with $\odd(T'_{i-1}) = 0$.

  We complete the proof by following the constructions of $\rho(F_i, 1)$ and $\rho(F_i, 0)$ as introduced in Section~\ref{sec:tplr}.
  First, note that
  \[F_i = t_i\ \otimes\ (X_i\ \oplus\ Y_i\ \oplus\ Z_i)\]
  has the universal cut vertex $t_i$ and that $F_i - t_i$ consists of the three mutually isomorphic connected components
  \begin{align*}
    X_i &= x_i\ \otimes\ (F_{i-1}\ \oplus\ u_i),\\
    Y_i &= y_i\ \otimes\ (F_{i-1}\ \oplus\ v_i),\\
    Z_i &= z_i\ \otimes\ (F_{i-1}\ \oplus\ w_i).
  \end{align*}
  Each of them has a universal cut vertex $x_i$, $y_i$ or $z_i$ and all of $X_i - x_i$, $Y_i - y_i, Z_i - z_i$ consist of the non-trivial connected component $F_{i-1}$ and an isolated vertex $u_i$, $v_i$ or $w_i$.
  This means that the leaf root construction for $X_i, Y_i, Z_i$ works according to the case $s=1$ and, up to isomorphism, we therefore have
  \begin{align*}
    (T, k) &= \rho(X_i, 1) = \rho(Y_i, 1) = \rho(Z_i, 1) \text{ with}\\
    k &= k_{i-1} + 2(1-\odd(T_{i-1})) = 2^{i+1}+1 \text{ and}\\
    (T', k') &= \rho(X_i, 0) = \rho(Y_i, 0) = \rho(Z_i, 0) \text{ with}\\
    k' &= k'_{i-1} + 2(1-\odd(T'_{i-1})) = 2^{i+1} + 2^{i-1}.
  \end{align*}
  Like in the proof of Proposition~\ref{thm:correct}, we get that a diametral path has length $\ell_1$ in both constructions and, thus,
  \begin{align*}
    \diam(T) &= 2k - 3 + \odd(T_{i-1}) = 2^{i+2}-1 \text{ and}\\
    \diam(T') &= 2k' - 3 + \odd(T'_{i-1}) = 2^{i+2} + 2^{i} - 3.
  \end{align*}
  Hence, $\odd(T) = \odd(T') = 1$.

  Now get back to the leaf root construction for $F_i$.
  If $p=1$ then our procedure handles the only case with $s \geq 2$. 
  We get
  \begin{align*}
    (T_i, k_i) &= \rho(F_i, 1) \text{ with}\\
    k_i &= 2k - 1 - 2\odd(T)^2 = 2^{i+2} - 1
  \end{align*}
  as claimed and, like in the proof of Proposition~\ref{thm:correct}, a diametral path is of length $\ell_{1,2}$, thus,
  \begin{align*}
    \diam(T_i) &= k_i + 2k - 2\odd(T) - 3\\
               &= (2^{i+2} - 1) + 2(2^{i+1}+1) - 5\\
               &= 2^{i+3} - 4,
  \end{align*}
  which means $\odd(T_i) = 0$.
  If $p = 0$ then the procedure works the case with $s > 2$ and where $k-\odd(T)$ equally stands for all three (isomorphic) $k$-leaf roots $T$ of the sub graphs $X_i$, $Y_i$, $Z_i$.
  This means that
  \begin{align*}
    (T'_i, k'_i) &= \rho(F_i, 0) \text{ with}\\
    k'_i &= 2k' - 2\odd(T')^2 = (2^{i+2} + 2^{i}) - 2\\
         &= (2^{i+2} - 1) + 2^{i} - 1 = k_i + 2^{i} - 1,
  \end{align*}
  as claimed.
  Moreover, like in the proof of Proposition~\ref{thm:correct}, a diametral path in $T'_i$ has length $\ell_{2,3}$.
  Here, we get
  \begin{align*}
    \diam(T'_i) &= k'_i - \odd(k'_i) + 2k' - 2\odd(T') - 2\\
               &= (2^{i+2} + 2^{i} - 2) + 2(2^{i+1} + 2^{i-1}) - 4\\
               &= 2^{i+3} + 2^{i+1} - 6,
  \end{align*}
  thus, $\odd(T_i) = 0$.
  This completes the proof.
\qed\end{proof}

\bibliographystyle{splncs04} 
\bibliography{leafpowers}
\end{document}